\documentclass[final]{autart}
\usepackage{amsmath}
\usepackage{amssymb}
\usepackage{makeidx}
\usepackage{multicol}
\usepackage[bottom]{footmisc}
\usepackage{graphicx}
\usepackage{color}
\usepackage{natbib}
\usepackage{epstopdf}
\usepackage{algorithm}
\usepackage{algcompatible}
\usepackage{calc}

\newtheorem{approxi}{Approximation}
\newtheorem{assump}{Assumption}

\newcommand{\app}{\begin{approxi}}
\newcommand{\eapp}{\end{approxi}}
\newcommand{\ass}{\begin{assump}}
\newcommand{\eass}{\end{assump}}
\newcommand{\teo}{\begin{thm}}
\newcommand{\eteo}{\end{thm}}
\newcommand{\corr}{\begin{cor}}
\newcommand{\ecorr}{\end{cor}}
\newcommand{\pro}{\begin{prop}}
\newcommand{\epro}{\end{prop}}
\newcommand{\lemma}{\begin{lem}}
\newcommand{\elemma}{\end{lem}}
\newcommand{\pb}{\begin{prob}}
\newcommand{\epb}{\end{prob}}
\newcommand{\df}{\begin{defn}}
\newcommand{\edf}{\end{defn}}
\newcommand{\rema}{\begin{rem}}
\newcommand{\erema}{\end{rem}}

\newcommand{\al}[1]{\begin{align} #1 \end{align}}
\newcommand{\nn}{\nonumber}

\newcommand{\tr}{\mathop{\rm tr}}  

\newcommand{\Ac}{ \mathcal{A}}

\newcommand{\Cc}{ \mathcal{C}}
\newcommand{\Dc}{ \mathcal{D}}
\newcommand{\Ec}{ \mathcal{E}}

\newcommand{\Hc}{ \mathcal{H}}
\newcommand{\Ic}{ \mathcal{I}}

\newcommand{\Pc}{ \mathcal{P}}

\newcommand{\Sc}{ \mathcal{S}}

\newcommand{\Es}{ \mathbb{E}}

\newcommand{\Ns}{ \mathbb{N}}

\newcommand{\Rs}{ \mathbb{R}}

\newcommand{\Zs}{ \mathbb{Z}}

\newcommand{\yv}{\mathrm{y}}
\newcommand{\uv}{\mathrm{u}}
\newcommand{\aaa}{\mathrm{a}}
\newcommand{\bbb}{\mathrm{b}}
\newcommand{\ccc}{\mathrm{c}}
\newcommand{\ev}{\mathrm{e}}

\newcommand{\pp}{\mathbf{p}}

\begin{document}

\begin{frontmatter}
\title{Nonparametric Identification of Kronecker Networks
}


\author[Padova]{Mattia Zorzi}\ead{zorzimat@dei.unipd.it}

\address[Padova]{Dipartimento di Ingegneria dell'Informazione, Universit\`a degli studi di
Padova, via Gradenigo 6/B, 35131 Padova, Italy}


\begin{abstract} We address the problem to estimate a dynamic network whose edges describe Granger causality relations and whose topology has a Kronecker structure. Such a structure arises in many real networks and allows to understand the organization of complex networks. We proposed a kernel-based PEM method to learn such networks. Numerical examples show the effectiveness of the proposed method. 
\end{abstract}
\end{frontmatter}

\section{Introduction}
In many applications we collect high dimensional data which describe systems having a wide amount of variables. Thus, there is the need to understand how those variables are ``connected each other'', i.e. we also need to learn from data the topology of the graphical model \citep{LAURITZEN_1996} representing  those interactions. It is worth noting that a graphical model is meaningful only if it is not complete, i.e. it has some degree of sparsity in terms of edges. It is clear that the network identification paradigm depends on the definition of ``interactions''. 
The latter can describe conditional dependence relations between the stochastic processes characterizing the variables of the system and whose resulting network is undirected   \citep{LINDQUIST_TAC13,LATENTG,8976281,RECS,zorzi2019graphical}. In the case of directed graphical models there is an edge from node $h$ to node $j$ if we need the process at node $h$ in order to reconstruct the process at node $j$. 
It is worth noting that such reconstruction can be both acausal \citep{MATERASSI1,MATERASSI2,v2020topology,sepehr2019blind,doddi2020estimating} and causal \citep{CHIUSO_PILLONETTO_SPARSE_2012,BSL_CDC}. In the latter case, we say that edges describes Granger causality relations \citep{GRANGER_CAUSALITY}. Finally, it is worth noting that a closely related topic is the problem to learn modules of a network only using data from neighbors nodes, see \cite{HOF1,HOF2}.

In this we paper we focus our attention on Kronecker networks \citep{leskovec2007scalable} whose topology  
possesses some properties emerging in many real networks, e.g. heavy-tailed degree distribution, small diameter, see \cite{leskovec2010kronecker} for more details. Such a structure is also useful to understand the organization of complex networks \citep{leskovec2009networks}. \cite{tsiligkaridis2013convergence} consider a network identification problem whose Kronecker structure is not only imposed on the topology of the network but also on the parameters of the model. The latter constraint, however, could be restrictive in some applications, e.g. for spatio-temporal MEG/EEG 
modeling \citep{bijma2005spatiotemporal}. A possible way to overcome this restriction is to impose the structure ``sum of Kronecker products'' on the parameters of the model \citep{QUARKS,tsiligkaridis2013covariance} at the price that we loose a meaningful definition of interaction, i.e. the resulting network is fully connected. In \cite{KRON,CDC_KRON}, instead, it has been considered the problem to learn undirected dynamic networks having the Kronecker structure on the topology, but not on the parameters of the model. Such a framework, however, has not been extended to the case of directed dynamic networks yet.

The present paper considers the problem to learn directed dynamic networks whose edges characterize Granger causality relations and whose topology has a Kronecker structure.

A well established system identification paradigm is the so called prediction error method (PEM), see \cite{LJUNG_SYS_ID_1999,SODERSTROM_STOICA_1988}. Within this framework, candidate models are described in fixed parametric model structures, e.g. ARMAX. Thus, the main difficulty is the correct choice of the best model structure which is usually performed by AIC and BIC criteria \citep{AKAIKE_1974}. Kernel-based PEM methods aim to overcome this issue, see \cite{PILLONETTO_DENICOLAO2010,EST_TF_REVISITED_2012,CHIUSO201624,MU2018381,ljung2020shift,9143975}. More precisely, the candidate model, described through the predictor impulse responses, is searched in an infinite dimensional nonparametric model class; this is clearly an ill posed problem because we have a finite set of measured data. However, it can be made into a well posed one using a penalty term which favors models with specific features. In the Bayesian viewpoint,  this is equivalent to introduce an a priori Gaussian probability (prior) on the unknown model. Hence, the prior distribution is characterized by the covariance (i.e. kernel) function.

The contribution of this paper is the introduction of a kernel-based PEM method to estimate Kronecker networks. The a priori information is that the impulse responses must be Bounded Input Bounded Output (BIBO) stable and the topology of the network respects the Kronecker structure. We derive the corresponding kernel functions by using the maximum entropy principle.

The kernel functions are characterized by the decay rate of the predictor impulse responses and by the number of edges in the network. These features are not known and characterized by the so called hyperparameters. We estimate them by minimizing the negative log-marginal likelihood of the measured data, see \cite{RASMUSSEN_WILLIAMNS_2006}. Moreover, we show that the negative log-marginal likelihood automatically sets some hyperparameters in such a way that the a priori information is that the network has a sparse Kronecker structure. The degree of sparsity depends on the variance of the noise process affecting each node.

Finally, we provide some numerical experiments on synthetic data in order to test the effectiveness of the proposed method. Moreover, we use the latter for learning Granger causality relations in a bike sharing system.

The paper is organized as follows. In Section  \ref{section_kron} we introduce the Kronecker model corresponding to a network characterized in terms of Granger causality relations. In Section \ref{section_pb_formulation} we introduce the Kronecker network identification problem. 
In Section \ref{section_kernel} we derive the maximum entropy kernels inducing BIBO stability and a Kronecker structure for the topology of the network, while
Section \ref{sec:opt_procedure} is devoted to the estimation of the hyperparameters.
In Section \ref{section_simulation} we present some numerical examples. Finally, we draw the conclusions in Section  \ref{section_conclusions}.
\subsection*{Notation}
 $\Ns$ is the set of natural numbers. Given a finite set $\Ic$,
$| \Ic |$ denotes its cardinality.
$ \Es[\cdot]$ denotes the expectation, while $ \Es[\cdot|\cdot]$ denotes the conditional mean. Given three (possibly infinite dimensional) random vectors $\aaa$, $\bbb$ and $\ccc$ we say that $\aaa$ is conditionally independent of $\bbb$ given $\ccc$ if  $\Es[\aaa|\bbb,\ccc] =  \Es[\aaa|\ccc]$. Given $G\in\Rs^{n\times p}$, $[G]_{i,j}$ denotes the entry of $G$ in position $(i,j)$; $G>0$ ($G\geq 0$) means that $G$ is a positive (semi-)definite matrix. $\Dc_p$ denotes the vector space of diagonal matrices of dimension $p$. $\ell_2(\Ns)$ denotes the space of $\Rs$-valued infinite length sequences, which we think as infinite dimensional column vectors 
$g:=[g_1 \; g_2 \; \ldots \; g_j \; \ldots]^\top$,  $g_k\in\Rs$, $k\in\Ns$, such that
$\|g\|_2:=\sqrt{\sum_{k=1}^\infty |g_k|^2}<\infty$. 
$\ell_2^{p\times n}(\Ns)$ is the space of matrices of sequences in $\ell_2(\Ns)$
\al{\Phi=\left[\begin{array}{ccc} (\phi^{[11]})^\top & \ldots  & (\phi^{[1n]})^\top \\ \vdots  & \ddots & \vdots \\ (\phi^{[p1]})^\top & \ldots &  (\phi^{[pn]})^\top\end{array}\right]\nn}
where $\phi^{[ij]}\in\ell_2 (\Ns)$, $i=1\ldots p$ and $j=1\ldots n$. 
$\ell_1(\Ns)$ denotes the space of $\Rs$-valued infinite length sequences $g$  
such that $\|g\|_1:=\sum_{k=1}^\infty |g_k|<\infty$. $\ell_1^{p\times n}(\Ns)$ is defined in similar way. 
$\Sc_2(\Ns)$ denotes the space of symmetric infinite dimensional matrices $K$ such that  
$\| K\|_2:=\sqrt{\sum_{i,j=1}^\infty | [K]_{i,j}|^2}<\infty$. $\Sc_2^{p}(\Ns)$ is the space of symmetric infinite dimensional matrices 
\al{K=\left[\begin{array}{cccc} K^{[11]} & K^{[12]} & \ldots  & K^{[1p]} \\ 
K^{[12]} & K^{[22]} &  & K^{[2p]} \\ 
\vdots  & & \ddots & \vdots \\ K^{[1p]} & K^{[2p]} &\ldots & K^{[pp]}\end{array}\right]\nn}
where $K^{[ij]} \in\Sc_2(\Ns)$, $i,j=1\ldots p$.
Given $\Phi\in\ell_2^{p \times n}(\Ns)$, $\Psi\in\ell_2^{m \times n}(\Ns)$ and $K\in\Sc_2^{n}(\Ns)$, the products $\Phi \Psi^\top$ and $\Phi K \Psi^\top$ are understood as $p \times m$ matrices whose entries are limits of infinite sequences \citep{INFINITE_MATRICES}. Given $g\in \ell_2^{n\times 1}(\Ns)$ and $K\in \Sc_2^{n}(\Ns)$,  $\|g\|^2_{K^{-1}}:=g^\top K^{-1}g$. The definition is similar in the case that $g$ and $K$ have finite dimension. With some abuse of notation the symbol $z$ will denote both the complex variable as well as the shift operator 
$z^{-1} y(t) :=y(t-1)$. 
Given a stochastic process $y=\{y(t)\}_{t\in\Zs}$, with some abuse of notation, $y(t)$ will both denote a random vector and its sample value. 
From now on the time $t$ will denote \emph{present} and we shall talk about \emph{past} and \emph{future} with respect to time $t$. With this convention in mind,  \al{\yv^- =\left[\begin{array}{ccc}  y(t-1)^\top & y(t-2)^\top & \ldots  \end{array}\right]^\top\nn}
denotes the (infinite length) past data vector of $y$ at time $t$. We denote as $\mathrm{supp}(G(z))$ the support of the transfer matrix $G(z)$, i.e. the entries of $G(z)$ different from the null function correspond to entries equal to one in $\mathrm{supp}(G(z))$, otherwise the latter are equal to zero.   

\section{Kronecker Models} \label{section_kron}
We consider the following nonparametric model
\al{  \label{OEmodel_start}y(t)=& G(z)y(t)+F(z)u(t)+e(t)}
where $y(t)\in \Rs^{p_1p_2}$ is the output, $u(t)\in\Rs^m$ is the input, $e(t)$ is zero mean normalized white Gaussian noise (WGN) whose components have variance $\sigma_{hk}^2$ with $h=1\ldots p_1$ and $k=1\ldots p_2$; $G(z)=\sum_{t=1}^\infty G_t z^{-t}$ is a BIBO stable transfer matrix of dimension $p_1 p_2\times p_1 p_2$, while $F(z)=\sum_{t=1}^\infty F_t z^{-t}$ is  a BIBO stable transfer matrix of dimension $p_1p_2\times m$.  
The minimum  variance one-step ahead predictor of $y(t)$ based on the past data $\yv^-$, denoted by $\hat y(t|t-1)$, is given by 
\al{ \hat y(t|t-1)&= G(z)y(t)+F(z)u(t).\nn} Thus, $e(t)$ is the one-step-ahead prediction error 
$$
e(t)  = y(t) - \hat y(t|t-1).
$$
We assume that 
\al{\label{cond_supp}\begin{aligned}\mathrm{supp}(G(z))&=E_1 \otimes E_2\\
\mathrm{supp}(F(z))&=A_1  \otimes A_2 \end{aligned}}
where $E_1$, $E_2$, $A_1$ and $A_2$ denote the adjacency matrices of dimension $p_1\times p_1$, $p_2\times p_2$, $p_1\times 1$ and $p_2\times m$, respectively. Let $y_{hk}$ denote the scalar process of $y$ in position $(h-1)p_2+k$ with $h\in I_1:=\{1\ldots p_1\}$  
and $k\in I_2:=\{1\ldots p_2\}$. In a similar way, $\yv^-_{hk}$ denotes the past of $y_{hk}$ at time $t$. Finally, $u_i$ denotes the scalar component of $u$ in position $i\in I_u:=\{1\ldots m\}$ and $\uv_i^-$ the $i$-th component of $\uv^-$.
Throughout the paper we assume that $E_1$, $E_2$, $A_1$ and $A_2$ are different from the null matrix.

It is possible to describe the structure of model  (\ref{OEmodel_start})-(\ref{cond_supp}) using a Bayesian network 
\citep{LAURITZEN_1996}. The nodes correspond to the scalar processes $y_{hk}$'s, with $h\in I_1$ and $k\in I_2$, while $u_i$'s, with $i\in I_u$, represent exogenous variables affecting some nodes of the network. Then, the connections among the nodes obey the following rules: 
\begin{itemize}
\item there is a  directed link from  node $y_{jl}$ to  node  $y_{hk}$  if
\al{
 \Es[y_{hk}(t) | \yv^-,\uv^-]  \neq  \Es[y_{hk}(t) |\yv_{\tilde j\tilde l}^-,\tilde j\neq j, \tilde l\neq l,\uv^-],\nonumber
}
 i.e. if $\yv_{jl}^-$ is needed to predict $y_{hk}(t)$  given $\uv^-$ and $\yv_{\tilde j\tilde l}^-$, for any  $\tilde j\in I_1\setminus \{j\}$ and $\tilde l\in I_2\setminus\{l\}$; In this case, we shall say $y_{jl}$ conditionally Granger causes $y_{hk}$. 
 \item there is a  directed link from $u_{i}$ to  node  $y_{hk}$ if
\al{
 \Es[y_{hk}(t) | \yv^-,\uv^-]  \neq  \Es[y_{hk}(t) |\yv^-,\uv^-_{\tilde i}, \tilde i\neq i]\nonumber
}
 i.e. if $\uv_{i}^-$ is needed to predict $y_{hk}(t)$  given $\yv^-$ and $\uv_{\tilde i}^-$, for any  $\tilde i\in I_u\setminus \{i\}$; In this case, we shall say $u_{i}$ conditionally Granger causes $y_{hk}$. 
 \end{itemize} In such a network we can recognize $p_1$ modules containing $p_2$ nodes and sharing the same topology described by $E_2$, while the interaction among these $p_1$ modules is described by $E_1$. Finally, $A_1$ describes the topology among  the input and the modules; $A_2$ describes the topology among  the input components and the nodes in a module which is the same in each module. These facts are formalized by the next proposition.

\pro\label{propo_cond_gr} Consider model (\ref{OEmodel_start})-(\ref{cond_supp}). Let $y^\star_h$ be the process obtained by stacking $y_{hk}$ with $k\in I_2$ and $y^\dag_k$ be the process obtained by stacking $h\in I_1$. Then:
\begin{itemize}
\item $y_j^\star$ conditionally Granger causes $y^\star_h$ if and only if $[E_1]_{h,j}=1$;
\item  $y_l^\dag$ conditionally Granger causes $y^\dag_k$ if and only if $[E_2]_{k,l}=1$;
\item $u$ conditionally Granger causes $y^\star_h$ if and only if $[A_1]_{h}=1$;
\item $u_i$ conditionally Granger causes $y^\dag_k$ if and only if $[A_2]_{k,i}=1$.
\end{itemize}
\epro
\begin{proof}  We only prove the first claim, the proof of the others is similar. Let $g^{[hk,jl]}\in\ell_1(\Ns)$ denote the impulse response in $G(z)$ having input $y_{jl}$ and output $y_{hk}$.  By (\ref{cond_supp}), we have 
 \al{\label{supp_1}g^{[hk,jl]}\neq 0, \, \hbox{ for some } k,l\in I_2\iff [E_1]_{h,j}=1.}
 From (\ref{OEmodel_start}) we have that 
 \al{y_{hk}&(t)=\sum_{\substack{j\in I_1\\ l\in I_2}}
 [G(z)]_{(h-1)p_2+k,(j-1)p_2+l}y_{jl}(t) \nn\\
 &+\sum_{i\in I_u}  [F(z)]_{(h-1)p_2+k,i}u_i(t) + e_{hk}(t).\nn} Thus, condition (\ref{supp_1}) is equivalent to 
 \al{\label{cond_prop1}\Es[y_{hk}(t) | \yv^-,\uv^-]  \neq  \Es[y_{hk}(t) |\yv_{\tilde j\tilde l}^-,\tilde j\neq j, \tilde l\neq l,\uv^-],\nn\\\hbox{for some } k,l\in I_2 \iff [E_1]_{h,j}=1.  } By (\ref{cond_prop1}), we have \al{\label{supp_2}\Es[y_{hk}(t)|\uv^-,\yv^-]\neq \Es[y_{hk}(t)|\yv^{\star -}_{\tilde j}, \tilde j\neq j,\uv^-],  \nn\\ \hbox{for some } k\in I_2 \iff [E_1]_{h,j}=1. } Staking condition (\ref{supp_2}) for any $k\in I_2$, we obtain
 \al{&\Es[y_{h}^\star(t)|\uv^-,\yv^-]\neq \Es[y_{h}^\star(t)|\uv^-,\yv^{\star -}_{\tilde j}, \tilde j\neq j] \nn\\ &\hspace{4cm} \iff [E_1]_{h,j}=1\nn} which proves the claim.
\qed
\end{proof}

\begin{figure}
\includegraphics{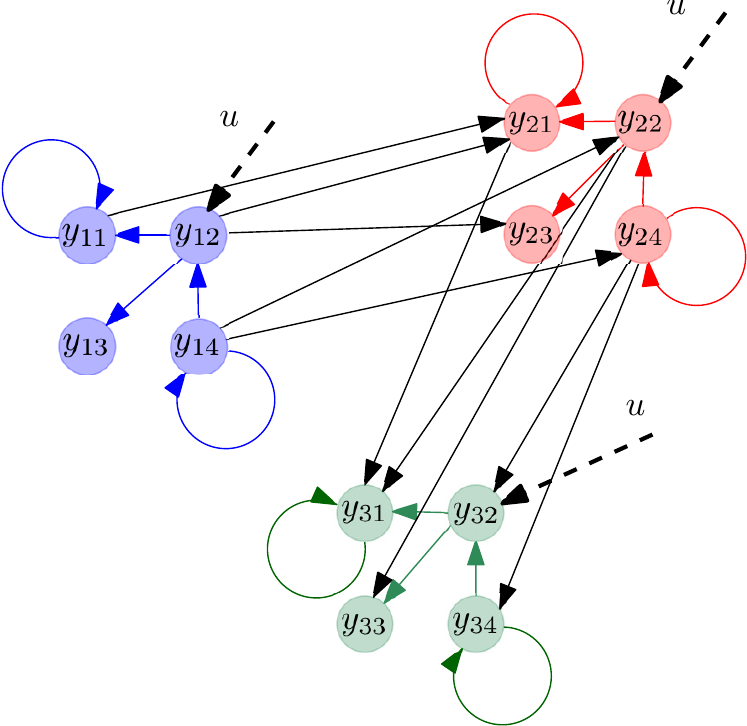}
\caption{A example of dynamic Kronecker network with $p_1=3$ modules (red, blue and green) composed by $p_2=4$ nodes and $m=1$ exogenous component.}\label{graph_ex}
\end{figure}
It is worth noting that condition (\ref{cond_supp}) is weaker than $G(z)$ and $F(z)$ admit a Kronecker decomposition, e.g. $G(z)=G_1(z)\otimes G_2(z)$. Therefore, we do not constrain the dynamics in each module to be the same. An example of Kronecker network is depicted in Figure \ref{graph_ex}. Here, we have 
$$y=[\,  y_{11} \,y_{12} \,y_{13} \,y_{14} \,y_{21} \,y_{22} \,y_{23} \,y_{24} \,y_{31} \,y_{32} \,y_{33} \,y_{34} \,]^\top$$ and the adjacency matrices are: 
{\small \al{E_1&=\left[\begin{array}{ccc}0 & 0 & 0 \\1 & 0 & 0 \\0 & 1 & 0\end{array}\right], \, E_2=\left[\begin{array}{cccc}1 & 1 & 0 & 0 \\0  & 0 & 0 & 1  \\0 & 1 & 0 &0 \\ 0 & 0 & 0 &1 \end{array}\right],\, A_1=\left[\begin{array}{c}1  \\1 \\1 \end{array}\right],\, A_2=\left[\begin{array}{c}0  \\1 \\0  \\ 0\end{array}\right].\nn}}

\subsection{Motivating examples}\label{sec_examples}
{\em Dynamic spatio-temporal modeling.} Consider a nonparametric time-varying model of the form
\al{ \label{mod_LTV}x(s)=\tilde G_s(q) x(s)+\tilde F_s(q) w(s)+v(s),\; \; s\in\Zs} 
where $x(s)\in\Rs^{p2}$ is the output, $w(s)\in\Rs^{\tilde m}$ is the input, $v(t)$ is normalized WGN and $\{\tilde G_s(q),\, s\in\Zs\}$, $\{\tilde F_s(q),\, s\in\Zs\}$ are sequences of BIBO stable (strictly causal) transfer matrices of dimension $p_2\times p_2$ and $p_2\times \tilde m$, respectively, such that 
\al{\tilde G_s(q)= \tilde G_{s+p_1}(q), \; \; \tilde F_s(q)=\tilde F_{s+p_1}(q)}
where $q^{-1}$ is the shift operator corresponding to time variable $s$, e.g. $q^{-1}x(s)=x(s-1)$. Then, if we define 
\al{ y(t)&=[\, x((t-1)p_1+1)^\top\, \ldots\; x(tp_1)^\top\,]^\top\in\Rs^{p_1p_2}\nn\\
u(t)&=[\, w((t-1)p_1+1)^\top\, \ldots\; w(tp_1)^\top\,]^\top\in\Rs^{p_1\tilde m}\nn\\
e(t)&=[\, v((t-1)p_1+1)^\top\, \ldots\; v(tp_1)^\top\,]^\top\in\Rs^{p_1p_2},\nn}
we can write (\ref{mod_LTV}) as (\ref{OEmodel_start}) where $G(z)$ and $F(z)$ are BIBO stable transfer matrices of dimension $p_1p_2\times p_1p_2$ and $p_1p_2\times p_1\tilde m$, respectively. Moreover, 
{\scriptsize \al{G(z)=\left[\begin{array}{ccccc}G_{1,1}(z) &   G_{1,2}(z) & \ldots   &\ldots &  G_{1,p_1}(z)  \\ zG_{2,1}(z) & G_{2,2}(z) & \ddots & & \vdots  \\
\vdots  & \ddots  & \ddots & \ddots &\vdots \\
\vdots & & \ddots & \ddots &   G_{p_1-1,p_1}(z) \\
zG_{p_1,1}(z) & \ldots & \ldots & zG_{p_1,p_1-1}(z) & G_{p_1,p_1}(z)\end{array}\right]\nn}}where $G_{hj}(z)=\sum_{t=1}^\infty G_{hj,t}z^{-t}$ is a transfer matrix of dimension $p_2\times p_2$ and $z^{-1}$ is the shift operator corresponding to time variable $t$; $F(z) $ is defined in a similar way: it is  composed by the transfer matrices $F_{hj}(z)=\sum_{t=1}^\infty F_{hj,t}z^{-t}$, $h,j=1\ldots p_1$, of dimension $p_2\times \tilde m$ and the blocks in the strictly lower block triangular part are multiplied by $z$. In plain words, we have written (\ref{mod_LTV}) as a time-invariant model through $y(t)$ and $u(t)$. It is interesting to note that the time variables $s$ and $t$ play a different role; indeed, there is a decimation relationship between them and the decimation factor is $p_1$. If $G(z)$ and $F(z)$ satisfy condition (\ref{cond_supp}), then by Proposition \ref{propo_cond_gr} we have that $E_2$ is the matrix describing the conditional Granger causality relations among the components of $x(s)$, while $A_2$ is the matrix describing the conditional Granger causality relations among the components of $w(s)$ and $x(s)$. Moreover, $E_1$ describes the conditional Granger causality relations of $x(s)$ over the period $1\ldots p_1$, while $A_1$ describes the conditional Granger causality relations among $x(s)$ and $w(s)$ over the period $1\ldots p_1$. We conclude that (\ref{OEmodel_start}) can be understood as a dynamic spatio-temporal model. In particular, if we take 
\al{\label{cond_special} G(z)=\bar G, \;\; F(z)= \bar F,}
i.e. constant matrices, then (\ref{OEmodel_start}) describes a spatio-temporal model without dynamic. The latter has been used to model magnetoencephalography (MEG) measurements for mapping brain activity    
\citep{bijma2005spatiotemporal}: $z(k,t)=[\,x_k((t-1)p_1+1) \, \ldots \,x_k(tp_1)\,]^\top$ models the measurements corresponding to the $k$-th
 brain region during the $t$-th trial of length $p_1$, while $s(k,t)=[\,w_k((t-1)p_1+1) \, \ldots \,w_k(tp_1)\,]^\top$ is the stimulus applied to the $k$-th region during the $t$-th trial. These trials are assumed independent. Moreover, they are identically distributed because in each trial the task performed by the patient is the same and the way that the patient reacts to the stimulus is assumed to be the same. In our framework we do not impose the restrictive condition (\ref{cond_special}) and thus we let dependence among the trials. This means that the proposed model takes into account the fact that the trials are taken in a sequential way, and thus the fact that the behaviors of the patient in the next trial can be influenced by the previous one. Finally, such a dynamic spatio-temporal model could be also used for modeling the conditional  Granger causality relations among the number rental bikes in a bikesharing system and the corresponding weather information, see  Section \ref{sec_bike}.

{\em Multi-task learning.} Consider a network composed by $p_2\in\Ns$ agents (i.e. each node correspond to an agent). One agent is described by a model having an input and output. Our aim is to model such a network under $p_1\in\Ns$ heterogeneous conditions (i.e. tasks). The $k$-th agent under the $h$-th task 
has input $u(t)$ and output $y_{hk(t)}$. Then, we could model the network through the $p_1$ models:
\al{\label{eq_multitask}y_h^\star(t)=G_h(z) y_h^\star+ F_h(z)u(t)+e_h(t), \; \; h\in I_1 }
where $y_h^\star(t)$ has been defined in Proposition \ref{propo_cond_gr}. A more flexible approach is to model all the models in (\ref{eq_multitask}) together in order to exploit commonalities and differences across the tasks, see \cite{yu2009large}. More precisely, we consider the model
\al{y(t)=G(z)y(t)+F(z)u(t)\nn}
where $y(t)$ is obtained by stacking $y_h^\star(t)$ with $h\in I_1$. 
Under the assumption that (\ref{cond_supp}) holds, and in view of Proposition \ref{propo_cond_gr}, we have that $E_1$ describes the conditional Granger causality relations among the tasks;  $E_2$ describes the conditional Granger causality relations among the agents in each task; $A_1$ describes the conditional Granger causality relations from the input to the tasks;  $A_2$ describes the conditional Granger causality relations from the input components to the agents in each task.

 \section{Problem Formulation}\label{section_pb_formulation}
Assume to measure the data $\{y(t),u(t)\}_{t=0,..,N}$ generated by (\ref{OEmodel_start}). In this section, we address the problem of  estimating $G(z)$ and $F(z)$ which admit the Kronecker product decomposition in (\ref{cond_supp}).

The transfer matrices  $G(z)$, $F(z)$ are parametrized in terms of their impulse response coefficients $G_t$ and $F_t$. In particular, defining 
 $g^{[hk,jl]}\in \ell_1(\Ns)$ to be the impulse response from input $(j-1)p_2+l$ to output $(h-1)p_2+k$,  we have:  
  \al{   g^{[hk,jl]}:= &\left[
                  \begin{array}{cc}
                    [G_1]_{(h-1)p_2+k,(j+1)p_2+l}   & \ldots \end{array} \right. \nn\\ 
                   &\hspace{0.2cm} \left. \begin{array}{ccc}  \ldots & [G_s]_{(h-1)p_2+k,(j+1)p_2+l}&\ldots \\
                  \end{array}
                \right]^\top ; \nn}
defining  $f^{[hk,i]}\in \ell_1(\Ns)$ to be the impulse response from input $u_i$ to output $(h-1)p_2+k$, we have
  \al{  f^{[hk,i]}:= \left[
                  \begin{array}{cccc}
                    [F_1]_{(h-1)p_2+k,i}  & \ldots & [F_s]_{(h-1)p_2+k,i,i}&\ldots \\
                  \end{array}
                \right]^\top . \nn}                
The coefficient vectors $\theta_g^\top\in\ell_1^{1\times p_1^2p_2^2}(\Ns)$ and $\theta_f^\top\in\ell_1^{1\times p_1p_2m}(\Ns)$ are defined  as follows:
\al{  \theta_g^\top =& \left[
                   \begin{array}{ccc|c}
                    (\theta_g^{[11]})^\top & \ldots  & (\theta_g^{[1p_2]})^\top & \ldots\\
                   \end{array}
                 \right.\nn\\
                 & \hspace{0.2cm} \left.
                   \begin{array}{c|ccc}
                      \ldots & (\theta_g^{[p_11]})^\top & \ldots &(\theta_g^{[p_1p_2]})^\top \\
                   \end{array}
                 \right]\nn}
\al{   \theta_f^\top =& \left[
                   \begin{array}{ccc|c}
                    (\theta_f^{[11]})^\top & \ldots  & (\theta_f^{[1p_2]})^\top & \ldots\\
                   \end{array}
                 \right.\nn\\
                 & \hspace{0.2cm} \left.
                   \begin{array}{c|ccc}
                      \ldots & (\theta_f^{[p_11]})^\top & \ldots &(\theta_f^{[p_1p_2]})^\top \\
                   \end{array}
                 \right]\nn}
where 
\al{ \theta_g^{[hk]^\top} =& \left[
                   \begin{array}{ccc|c}
                    (g^{[hk,11]})^\top & \ldots  & (g^{[hk,1p_2]})^\top & \ldots\\
                   \end{array}
                 \right.\nn\\
                 & \hspace{0.2cm} \left.
                   \begin{array}{c|ccc}
                      \ldots & (g^{[kh,p_11]})^\top & \ldots &(g^{[hk,p_1p_2]})^\top \\
                   \end{array}
                 \right]\nonumber\\
                 \theta_f^{[hk]^\top} =& \left[
                   \begin{array}{ccc}
                    (f^{[hk,1]})^\top & \ldots  & (f^{[hk,m]})^\top \\
                   \end{array}
                 \right]\nonumber.
                 }   
The measured data $y(1)\ldots y(N)$ are stacked in the vector $\yv^+$ as follows
  \al{  \yv^+ =&\left[
          \begin{array}{ccc|c|ccc}
            \yv_{11}^{+\top} &  \ldots & \yv_{1p_2}^{+\top}  & \ldots     & \yv_{p_11}^{+\top}&  \ldots & \yv_{p_1p_2}^{+\top} \\
          \end{array}
        \right]^\top\nn}
        where 
         \al{\yv^+_{hk} =&\left[
          \begin{array}{ccc}
            y_{hk}(N)^\top &  \ldots & y_{hk}(1)^\top   \\
          \end{array}
        \right].\nn
     }
         The vector $\ev^+$ is defined analogously.
Let us also introduce the {\em Toeplitz} matrices  $\phi_{hk} \in\ell_2^{N\times 1}(\Ns)$, $\psi_{i} \in\ell_2^{N\times 1}(\Ns)$:
 \al{ &[\phi_{hk}]_{sn}:=y_{hk}(N-s- n+1)\nn\\
 &[\psi_{i}]_{sn}:=u_{i}(N-s- n+1)
 \nn}
with $s=1\ldots N$ and $n\in\Ns$. Then, we define the  regression matrices   $\Phi\in\ell_2^{p_1p_2N \times p_1^2p_2^2} (\Ns)$ and $\Psi\in\ell_2^{p_1p_2N \times p_1p_2m} (\Ns)$  as:
 \al{ \Phi&=I_{p_1p_2} \otimes \left[
                       \begin{array}{ccccccc}
                         \phi_{11} & \ldots & \phi_{1p_2} & \ldots  & \phi_{p_1 1} & \ldots & \phi_{p_1p_2}\\
                       \end{array}
                     \right]\nn\\
                     \Psi&=I_{p_1p_2} \otimes \left[
                       \begin{array}{ccccccc}
                         \psi_{1} & \ldots & \psi_{m}\\
                       \end{array}
                     \right]\nn}   
so that, from (\ref{OEmodel_start}) the vector $\yv^+$ containing the measured output data satisfies the linear regression model \al{\label{LM} \yv^+= \Phi\theta_g+\Psi\theta_f+\ev^+ }
and $\hat\yv^+:= \Phi \theta_g+\Psi\theta_f$ is the one-step ahead predictor of $\yv^+$. It is worth noting that to construct $\Phi$ and $\Psi$  we need the remote past of the output and the input which is not available. Thus, model \eqref{LM} has to be approximated truncating $\Phi$ and $\Psi$ (and thus $\theta_g$ and $\theta_f$). This is equivalent to impose zero initial conditions, which is a reasonable approximation given the  decay, as a function of $t$, of the coefficients $G_t$, $F_t$ by the BIBO stability conditions. 

\rema \label{remark_spatio}In the case that we consider the dynamic spatio-temporal model of Section \ref{sec_examples}, the regression matrices in (\ref{LM}) are different because the strictly lower block triangular parts of $G(z)$ and $F(z)$ do not contain a delay. More precisely, $\Phi$ and $\Psi$ are $p_1\times p_1$ block diagonal matrices where the block in position $(h,h)$, with $h\in I_1$, is, respectively,
\al{&I_{p_2}\otimes [\,\tilde\phi_{11} \,\ldots\, \tilde\phi_{h-1,p_2}\,  \phi_{h,1}\, \ldots \phi_{p_1,p_2}  \,]\nn\\
&I_{p_2}\otimes [\,\tilde \psi_{1} \,\ldots\, \tilde \psi_{(h-1)\tilde m}\,  \psi_{(h-1)\tilde m+1}\, \ldots  \psi_{p_1\tilde m}  \,]\nn} and the {\em Toeplitz} matrices  $\tilde \phi_{hk} \in\ell_2^{N\times 1}(\Ns)$, $\tilde \psi_{i} \in\ell_2^{N\times 1}(\Ns)$ are:
 \al{ &[\tilde \phi_{hk}]_{sn}:=y_{hk}(N-s- n+2)\nn\\
 &[\tilde \psi_{i}]_{sn}:=u_{i}(N- s- n+2).
 \nn}
\erema 

Therefore, our Kronecker identification problem can be formulated in terms of PEM as follows.
 \pb \label{problem} Assume to collect the measurements  $\{y(t),u(t) \}_{t=0\ldots N}$ and that the dimensions $p_1$ and $p_2$ are known. Find $\theta_g^\top\in\ell_1^{1\times p_1^2p_2^2}(\Ns)$ and $\theta_f^\top\in\ell_1^{1\times p_1p_2m}(\Ns)$ corresponding to a Kronecker model  minimizing the prediction error squared norm $\| \yv^+-\Phi \theta\|^2_{\Sigma^{-1}\otimes I_N}$ with $\Sigma=\mathrm{diag}( \sigma_{11}^2 \ldots\sigma_{1p_2}^2 \ldots\sigma_{p_11}^2 \ldots\sigma_{p_1 p_2})$. \epb Following the nonparametric Gaussian regression approach in \cite{PILLONETTO_2011_PREDICTION_ERROR},  
we model $\theta_g$ and $\theta_f$ as zero-mean processes
with kernels $K_g\in\Sc_2^{p_1^2p_2^2}(\Ns)$ and $K_f\in\Sc_2^{p_1p_2m}(\Ns)$, respectively. These kernels may depend upon some tuning parameters, usually called hyperparameters and denoted with  $\zeta$ hereafter. As illustrated in Section \ref{section_kernel},  according to the maximum entropy principle $\theta_g$ and $\theta_f$ will be modeled as  Gaussian and independent.
In the following $\Hc_{K}$ denotes the reproducing Hilbert space  \citep{ARONSZAJN1950} of deterministic functions on $\Ns$, associated with $K$ and with norm denoted by $\|\cdot \|_{K^{-1}}$. We assume that the past data $\yv^-$ neither affects the {\em a priori} probability on $\theta_g$  and $\theta_f$ nor carries information on $\zeta$ and $\Sigma$ \citep{PILLONETTO_DENICOLAO2010}, that is 
 \al{ \label{approx_pb}\mathbf{p} & (\yv^+,\theta_g,\theta_f,\yv^-|\zeta,\Sigma)\nn\\ 
& =\mathbf{p}(\yv^+|\theta_g,\theta_f,\yv^-,\zeta,\Sigma) \mathbf{p}(\theta_g,\theta_f|\yv^-,\zeta,\Sigma)\mathbf{p}(\yv^-|\zeta,\Sigma)\nn\\
& \approx
\mathbf{p}(\yv^+|\theta_g,\theta_f,\yv^-,\zeta,\Sigma) \mathbf{p}(\theta_g,\theta_f|\zeta,\Sigma)\mathbf{p}(\yv^-).}
Let $\hat \theta_g= \Es[\theta_g|\yv^+,\zeta,\Sigma]$ and $\hat \theta_f= \Es[\theta_f|\yv^+,\zeta,\Sigma]$  be the  minimum variance estimator, respectively, of $\theta_g$ and $\theta_f$ given $\yv^+$, $\zeta$ and $\Sigma$. It is well known  that $\hat \theta_g$ and $\hat \theta_f$ are almost surely solution to the following {\em Tikhonov}-type variational problem 
\al{ \underset{\substack{\theta_g\in\Hc_{K_g}\\\theta_f\in\Hc_{K_f}}}{\arg\min}\|\yv^+-\Phi   \theta_g -\Psi   \theta_f\|^2_{\Sigma^{-1}\otimes I_N}
  +\| \theta_g\|^2_{K_g^{-1}}+\| \theta_f\|^2_{K_f^{-1}}.\nn}
Moreover, almost surely: 
\al{  \begin{aligned}\hat \theta_g &= K_g\Phi^\top (\Phi K_g\Phi^\top+\Psi K_f\Psi^\top+\Sigma \otimes I)^{-1}\yv^+\\
\hat \theta_f&= K_f\Psi^\top (\Phi K_g\Phi^\top+\Psi K_f\Psi^\top+\Sigma \otimes I)^{-1}\yv^+\end{aligned}.\nn}
In what follows, $\hat G(z)$ and $\hat F(z)$ denote the transfer matrices corresponding to $\hat \theta_g$ and $\hat \theta_f$, respectively. The main task now is to design the kernels $K_g$ and $K_f$ in such a way that  $\hat G(z)$ and $\hat F(z)$ are almost surely BIBO stable while favoring a sparse Kronecker decomposition of $\mathrm{sup}(\hat G(z))$ and $\mathrm{sup}(\hat F(z))$.

 \section{Maximum entropy priors} \label{section_kernel}
 
The probability law for the joint process $\theta=[\,\theta_g^\top\; \theta_f^\top]^\top$ under desired constraints can be obtained by the maximum entropy principle. 
Indeed, maximum entropy solutions rely on ``information'' arguments which essentially state that the maximum entropy distribution is the one satisfying the given constraints and containing the largest amount of ``uncertainty'' \citep{1456693}. In plain words, this principle guarantees that this solution  does not satisfy additional constraints that are undesired and unnecessary. We shall make the rather mild assumption that the process $\theta$ is zero-mean and absolutely continuous with respect to the Lebesgue measure. We will see that  the optimal solution (i.e. maximizing the differential entropy) is a Gaussian 
process where $\theta_g$ and $\theta_f$ are independent. Then, we will also characterize the corresponding kernel functions.

We start with the constraints on $\theta_g$ inducing BIBO stability on $\hat G(z)$.  Let $P\in\Sc_2(\Ns)$ be a strictly positive definite kernel (in the sense of Moore) such that $[P]_{t,t}\leq \kappa t^\alpha e^{-\beta t}$, $t\in\Ns$, with $\kappa,\beta>0$ and $\alpha\in\Rs$.  Let also $\vartheta$ be a zero-mean process which satisfies the moment  constraint
    \al{ \Es[\| \vartheta\|_{P^{-1}}^2]\leq c\nn} 
  where $c\geq 0$. Then, for any $\varepsilon>1$ there exists $\bar \kappa_\varepsilon>0$ such that   the covariance function (kernel) $K$ of  $\vartheta$ satisfies \cite[Proposition 4]{BSL}:
 \al{[K]_{t,t}\leq \bar \kappa_\varepsilon t^{\alpha+\varepsilon} e^{-\beta t},\;\; t\in\Ns.\nn} 
 It is not difficult to see that the condition above on $P$  is satisfied by the kernels usually employed in the identification of dynamical models (e.g. stable spline, tuned/correlated and so on, see \cite{KERNEL_METHODS_2014}).   
Thus, we consider the constraints \al{ \label{constraint_p_s}  \Es[\| g^{[hk,jl]}\|^2_{P^{-1} }] \leq \tilde c_{hk,jl}, \; \forall\,h,j\in I_1,\;\forall\, k,l\in I_2}
where $c_{hk,jl} \geq 0$ and $P\in\Sc_2(\Ns)$ as above. Then, by (\ref{constraint_p_s}), the covariance of the $k$-th element of $g^{[hk,jl]}$ decays exponentially. Accordingly, its posterior mean is a BIBO stable transfer function. 

\rema The kernel $P$ induces a prior  about the decay rate of the impulse responses $g^{[hk,jl]}$. On the other hand, the hyperparameters of $P$ are estimated  from the measured data, see Section 5. Therefore, although $P$ could assign high probability to a small subset of BIBO stable systems, the latter is a  ``wise'' set because it has been extrapolated from the data.\erema

It remains to derive the constraints inducing a sparse Kronecker  decomposition as in (\ref{cond_supp}), i.e.:
\al{\label{cond_sp_1or}&[E_1]_{h,j}=0 \iff g^{[hk,jl]}=0, \; \forall\, k,l\in I_2\\
\label{cond_sp_2or}&[E_2]_{k,l}=0\iff g^{[hk,jl]}=0, \; \forall\, h,j\in I_1.} We consider the constraints:
\al{\label{constr1} \sum_{k,l\in I_2}\Es[\| g^{[hk,jl]}\|^2_{P^{-1} }] \leq c_{hj}^\prime, \; \forall\,h,j\in I_1;\\
\label{constr2} \sum_{h,j\in I_1}\Es[\| g^{[hk,jl]}\|^2_{P^{-1} }] \leq c_{kl}^{\prime\prime}, \; \forall\,k,l\in I_2 }
where $c_{hj}^\prime,c_{kl}^{\prime\prime} \geq 0$. First, notice that (\ref{constr1}) and (\ref{constr2}) imply (\ref{constraint_p_s}) provided that $\max\{c_{hj}^\prime,c_{kl}^{\prime\prime}\}\leq \tilde c_{hk,jl}$ for any $h,k,j,l$. Accordingly, the posterior mean of $G(z)$, under the constraints (\ref{constr1})-(\ref{constr2}), is BIBO stable. Let $c_{hj}^\prime=0$, by (\ref{constr1}) we have that $g^{[hk,jl]}$ $\forall\, k,l\in I_2$ are the null sequence in mean square and
 so are their posterior mean, i.e. the latter satisfy (\ref{cond_sp_1or}). 
Let $c_{kl}^{\prime \prime}=0$, by (\ref{constr2}) we have that $g^{[hk,jl]}$ $\forall\, h,j\in I_1$ are the null sequence in mean square and
 so are their posterior mean, i.e. the latter satisfy (\ref{cond_sp_2or}).
 In a similar way, condition (\ref{cond_supp}) on $F(z)$ can be written as:
 \al{ &[A_1]_{h}=0 \iff f^{[hk,i]}=0, \; \forall\, k\in I_2, \; \forall\, i\in I_u\nn\\
 &[A_2]_{k,i}=0\iff f^{[hk,i]}=0, \; \forall\, h\in I_1\nn.} Thus, the constraints inducing BIBO stability and a sparse Kronecker decomposition are:
\al{\label{constr1f} &\sum_{\substack{k\in I_2 \\i\in I_u}}\Es[\| f^{[hk,i]}\|^2_{R^{-1} }] \leq d_{h}^\prime, \; \forall\,h\in I_1;\\
&\label{constr2f} \sum_{h\in I_1}\Es[\| f^{[hk,i]}\|^2_{R^{-1} }] \leq d_{ki}^{\prime\prime}, \; \forall\,k\in I_2,\; \forall \, i\in I_u }
where $d_{h}^\prime,d_{ki}^{\prime\prime} \geq 0$ and $R\in\Sc_2(\Ns)$ is a strictly positive kernel such that $[R]_{t,t}\leq \kappa t^\alpha e^{-\beta t}$, $t\in\Ns$, $\kappa,\beta>0$ and $\alpha\in\Rs$.

In order to build the desired prior distribution we make use of the Kolmogorov extension Theorem, see  \cite{Oksendal}, and work with finite vectors extracted from process $\theta=[\, \theta_g^\top\;\theta_f^\top\,]^\top$. More precisely, consider a finite index set $\Ic=\Ic_g\times \Ic_f$ and $\Ic_g,\Ic_f\subset \Ns$. Let $\check \theta$ be the random vector whose components are extracted, from processes $\theta_g$ and $\theta_f$ according to the index sets $\Ic_g$ and $\Ic_f$, respectively. We denote by  $\pp_{\Ic}(\check \theta)$ the  probability density of $\check \theta$. In view of the Kolmogorov extension Theorem, process $\theta$ can be characterized by specifying the joint probability density $\pp_{\Ic}$ for all finite sets $\Ic \subset \Ns\times \Ns$. Thus, the maximum entropy process $\theta$ can be constructed building all the marginals   $\pp_{\Ic}$ using the maximum entropy principle, which can thus be extended by the Kolmogorov extension theorem.
 Such principle states that
among all the probability densities $\pp_{\Ic}$ satisfying the desired constraints, the optimal one should maximize the differential entropy \citep{COVER_THOMAS}.
Let $ P_{\Ic_g}$ and $R_{\Ic_f}$ be the kernel matrices whose entries are extracted from $P$ and $R$ according to $\Ic_g$ and $\Ic_f$, respectively. Therefore, the corresponding maximum entropy problem is 
 \al{ \label{ME_problem}\underset{\pp_{\Ic}\in\Pc}{\max}& \mathbf{H} (\pp_{\Ic})\nn\\
\hbox{s.t. } & \text{ (\ref{constr1}) (\ref{constr2}) (\ref{constr1f})  (\ref{constr2f}) hold for $\Ic$}}
where $\Pc$ is the class of probability densities in $\Rs^{|\Ic|}$
which are bounded and taking positive values and $\mathbf{H} (\pp_{\Ic})$ denotes the differential entropy of $\check \theta$. Moreover, the constraints in (\ref{ME_problem}) must be understood as follows: (\ref{constr1}), (\ref{constr2}), (\ref{constr1f}) and (\ref{constr2f}) hold with
$\theta$, $ P$ and $R$ replaced by $\check \theta$, $ P_{\Ic_g}$ and $R_{\Ic_f}$, respectively.

\teo\label{teo_ME} Under the assumption that $c_{hj}^\prime,c_{kl}^{\prime\prime},d_{h}^\prime,d_{ki}^{\prime\prime}>0$, $h,j\in I_1$,  $k,l\in I_2$ and $i\in I _u$, the unique optimal solution to the maximum entropy problem (\ref{ME_problem}) is 
such that 
$\check \theta=[\, \check \theta_g^\top \; \check \theta_f^\top\,]^\top$ is Gaussian with zero mean. Moreover, $\check \theta_g$ and $\check \theta_f$ are independent and with kernel matrix, respectively, 
\al{\label{optimal_K} \begin{aligned} \check K_g&= X_g\otimes P_{\Ic_g}\\ \check K_f&= X_f\otimes R_{\Ic_f}\\
X_g&=(\Lambda\otimes \Gamma)(\Lambda\otimes I_{p_2^2}+I_{p_1^2}\otimes \Gamma)^{-1}\\
X_f&=(\Pi\otimes \Omega)(\Pi\otimes I_{p_2m}+I_{p_1}\otimes \Omega)^{-1}\\
\Lambda &=\mathrm{diag}(\lambda_{11}\ldots \lambda_{1p_1} \ldots \lambda_{p_11} \ldots \lambda_{p_1 p_1})\\
\Gamma &=\mathrm{diag}(\gamma_{11}\ldots \gamma_{1p_2} \ldots \gamma_{p_21} \ldots \gamma_{p_2 p_2})\\
\Pi &=\mathrm{diag}(\pi_{1}\ldots \pi_{p_1}  )\\
\Omega &=\mathrm{diag}(\omega_{11}\ldots \omega_{1m} \ldots \omega_{p_21} \ldots \omega_{p_2 m}) \end{aligned}
}
 where $\lambda_{hj},\gamma_{kl}, \pi_h,\omega_{ki}>0$.
  \eteo
\begin{proof} We prove the claim by using the duality theory. The Lagrangian is 
\al{\label{Lagr}L(\pp_{\Ic},&\tilde \Lambda,\tilde \Gamma,\tilde \Pi,\tilde\Omega)=\mathbf{H}(\pp_{\Ic})\nn\\ &+\frac{1}{2}\sum_{h,j\in I_1}\tilde \lambda_{hj}\left(c_{hj}^\prime-
\sum_{k,l\in I_2}\Es[\| \check g^{[hk,jl]}\|^2_{P_{\Ic_g}^{-1} }] \right)\nn\\ &+\frac{1}{2}\sum_{k,l\in I_2}\tilde \gamma_{kl}\left( c_{kl}^{\prime\prime}-\sum_{h,j\in I_1}\Es[\| \check  g^{[hk,jl]}\|^2_{P_{\Ic_g}^{-1} }] \right)\nn\\ 
&+\frac{1}{2}\sum_{h\in I_1}\tilde \pi_{h}\left(d_{h}^\prime-
\sum_{\substack{k\in I_2\\ i \in  I_u}}\Es[\| \check  f^{[hk,i]}\|^2_{R_{\Ic_f}^{-1} }] \right)\nn\\ 
&+\frac{1}{2}\sum_{\substack{k\in I_2\\ i \in I_u}}\tilde \omega_{ki}\left( d_{ki}^{\prime\prime}-\sum_{h\in I_1}\Es[\| \check f^{[hk,i]}\|^2_{R_{\Ic_f}^{-1} }] \right)\nn\\
&= -\Es[\log( \pp_{\Ic})]+\frac{1}{2}\sum_{h,j\in I_1}\tilde \lambda_{hj} c_{hj}^\prime+\frac{1}{2}\sum_{k,l\in I_2} \tilde \gamma_{kl} c_{kl}^
{\prime\prime}\nn\\
&-\frac{1}{2}\sum_{\substack{h,j\in I_1\\k,l\in I_2}}(\tilde \lambda_{hj}+\tilde\gamma_{kl} )\Es[\| \check  g^{[hk,jl]}\|^2_{P_{\Ic_g}^{-1} }]\nn\\
&+\frac{1}{2}\sum_{h\in I_1}\tilde \pi_{h} d_{h}^\prime+\frac{1}{2}\sum_{\substack{k\in I_2\\ i\in I_u}} \tilde \omega_{ki} d_{ki}^{\prime\prime}\nn\\
&-\frac{1}{2}\sum_{\substack{h\in I_1,k\in I_2\\ i\in I_u}}(\tilde \pi_{h}+\tilde\omega_{ki} )\Es[\|\check  f^{[hk,i]}\|^2_{R_{\Ic_f}^{-1} }]\nn\\
&= -\Es[\log( \pp_{\Ic})]+\frac{1}{2}\sum_{h,j\in I_1}\tilde \lambda_{hj} c_{hj}^\prime+\frac{1}{2}\sum_{k,l\in I_2} \tilde \gamma_{kl} c_{kl}^{\prime\prime}\nn\\
&+\frac{1}{2}\sum_{h\in I_1}\tilde \pi_{h} d_{h}^\prime -\frac{1}{2}\Es[\|\check  \theta_g\|^2_{(\tilde \Lambda \otimes I+ I\otimes \tilde \Gamma )\otimes P_{\Ic_g}^{-1} }]\nn\\
&+\frac{1}{2}\sum_{\substack{k\in I_2\\ i\in I_u}} \tilde \omega_{ki} d_{ki}^{\prime\prime}-\frac{1}{2}\Es[\|\check  \theta_f\|^2_{(\tilde \Pi \otimes I+ I\otimes \tilde \Omega )\otimes R_{\Ic_f}^{-1} }]
}
where $\tilde \Lambda=\mathrm{diag}(\tilde \lambda_{11} \ldots \tilde \lambda_{1p_1} \ldots \tilde \lambda_{p_11 } \ldots \tilde \lambda_{p_1 p_1} )$, $\tilde \Gamma=\mathrm{diag}(\tilde \gamma_{11} \ldots\tilde \gamma_{1p_2} \ldots\tilde \gamma_{p_21} \ldots\tilde \gamma_{p_2 p_2} )$, $\tilde \Pi=\mathrm{diag}(\tilde  \pi_{1} \ldots  \ldots\tilde \pi_{p_1} )$, $\tilde \Omega=\mathrm{diag}(\tilde \gamma_{11} \ldots\tilde \gamma_{1m} \ldots\tilde \gamma_{p_21} \ldots\tilde \gamma_{p_2 m} )$; moreover, $\tilde \lambda_{h,j},\tilde\gamma_{kl}, \tilde \pi_h, \tilde \omega_{hi}\geq 0$ are the Lagrange multipliers corresponding to (\ref{constr1}), (\ref{constr2}), (\ref{constr1f}) and (\ref{constr2f}), respectively. It is not difficult to see that (\ref{Lagr}) is strictly concave in $\Pc$. Moreover, its point of maximum exists under  the assumptions that $\tilde \lambda_{hj}+\tilde \gamma_{kl}>0$ and $\tilde \pi_{h}+\tilde \omega_{ki}>0$, and it is Gaussian distributed: \al{\label{opt_pI}\pp_{\Ic}(\check\theta)= \pp_{\Ic_g}(\check\theta_g)\pp_{\Ic_f}(\check\theta_f)}
where 
\al{
\pp_{\Ic_f}(\check \theta_g)&=\frac{1}{\sqrt{(2\pi)^{p_1^2p_2^2|\Ic_g|} |\check K_g|}}\exp(-\frac{1}{2}\check\theta_g^\top \check K_g^{-1}\check\theta_g)\nn\\
\pp_{\Ic_f}(\check\theta_f)&=\frac{1}{\sqrt{(2\pi)^{p_1p_2m |\Ic_f|} |\check K_f|}}\exp(-\frac{1}{2}\check\theta_f^\top \check K_f^{-1}\check\theta_f)\nn}
and 
\al{\check K_g&=(\tilde \Lambda \otimes I+ I\otimes \tilde \Gamma )^{-1}\otimes P_{\Ic_g}\nn\\
\check K_f&=(\tilde \Pi\otimes I+ I\otimes \tilde \Omega )^{-1}\otimes R_{\Ic_f}.\nn} In the case that $\tilde \lambda_{hj}+\tilde \gamma_{kl}=0$ or $\tilde \pi_{h}+\tilde \omega_{ki}=0$ for some $h,k,j,l,i$, then the maximum cannot be attained in $\Pc$. Substituting the optimal form of $\pp_{\Ic}$ in the Lagrangian we obtain the dual function (up to terms not depending on $\tilde \Lambda$, $\tilde \Gamma$, $\tilde \Pi$ and $\tilde \Omega$):
\al{J(\tilde \Lambda,\tilde \Gamma, \tilde\Pi, \tilde\Omega )=J_f(\tilde \Lambda,\tilde \Gamma)+J_g(\tilde \Pi,\tilde \Omega)\nn}
where
\al{J_f(\tilde\Lambda,\tilde \Gamma)=-&\frac{|\Ic_g|}{2}\sum_{\substack{h,j\in I_1\\ k,l\in I_2}} \log(\tilde \lambda_{hj}+\tilde \gamma_{kl})\nn\\ &\hspace{0.5cm}+\frac{1}{2}\sum_{h,j\in I_1}\tilde \lambda_{hj} c_{hj}^\prime+\frac{1}{2}\sum_{k,l\in I_2} \tilde \gamma_{kl} c_{kl}^{\prime\prime}.
\nn\\
J_g(\tilde\Pi,\tilde \Omega)=-&\frac{|\Ic_f|}{2}\sum_{\substack{h\in I_1, \, k\in I_2\\ i\in I_u}} \log(\tilde \pi_{h}+\tilde \omega_{ki})\nn\\ &\hspace{0.5cm}+\frac{1}{2}\sum_{h\in I_1}\tilde \pi_{hj} d_{h}^\prime+\frac{1}{2}\sum_{\substack{k\in I_1\\ i\in I_u}} \tilde \omega_{ki} d_{ki}^{\prime\prime}.
\nn}  Hence, we can minimize $(\tilde \Lambda,\tilde\Gamma)$ and $(\tilde \Pi, \tilde \Omega)$ in an independent way. We start to minimize $J_f(\tilde \Lambda,\tilde \Gamma)$ over the open and unbounded set \al{\Cc=\{\, (\tilde \Lambda,\tilde \Gamma)\in \Dc  \hbox{ s.t. } \tilde \Lambda\geq0,\, \tilde \Gamma \geq0,\,  \tilde \Lambda \otimes I+I\otimes\tilde \Gamma>0 \, \}\nn} where $\Dc= \Dc_{p_1^2}\times \Dc_{p_2^2}$.

Next, we show that we can restrict the search of the minimum of $J_f$ over a compact set and thus by the Weierstrass' theorem, the dual problem admits solution. Take a sequence $(\tilde\Lambda^{(n)}, \tilde \Gamma^{(n)})\in \Cc$, $n\in\Ns$, such that $\|\tilde \Lambda^{(n)}\|\rightarrow \infty$ and/or $\|\tilde \Gamma^{(n)}\|\rightarrow \infty$ as $n\rightarrow \infty$. Then, $J_f(\tilde \Lambda^{(n)},\tilde \Gamma^{(n)})\rightarrow \infty$ because the linear terms dominate the logarithmic ones. Therefore, such a sequence cannot be an infimizing sequence. Accordingly, we can restrict the set $\Cc$ as
\al{\Cc_1=\{\, (\tilde \Lambda&,\tilde \Gamma)\in \Dc  \hbox{ s.t. } 0\leq  \tilde\Lambda\leq \lambda_{MAX}I,\nn\\ & 0\leq\Gamma \leq\gamma_{MAX}I,\, \tilde \Lambda\otimes I+I\otimes \tilde \Gamma>0   \, \}\nn}
where $\lambda_{MAX},\gamma_{MAX}>0$ are constants taken sufficiently large. Take a sequence $(\tilde\Lambda^{(n)}, \tilde \Gamma^{(n)})\in \Cc_1$, $n\in\Ns$, such that $\lambda^{(n)}_{hj}+\tilde \gamma_{kl}^{(n)}\rightarrow 0$ as $n\rightarrow \infty$ for some $h,j\in I_1$ and $k,l\in I_2$. Then, we have 
\al{J_f(\tilde \Lambda^{(n)},\tilde\Gamma^{(n)})\geq -\frac{|\Ic_g|}{2}\sum_{\substack{h,j\in I_1\\ k,l\in I_2}} \log(\tilde \lambda_{hj}+\tilde \gamma_{kl})\rightarrow \infty.\nn} Accordingly, it cannot be an infimizing sequence, and thus we can restrict the search to the set 
\al{\Cc_2=\{\, (\tilde \Lambda&,\tilde \Gamma)\in \Dc  \hbox{ s.t. } 0\leq  \tilde\Lambda\leq \lambda_{MAX}I,\nn\\ & 0\leq\Gamma \leq\gamma_{MAX}I,\, \tilde \Lambda\otimes I+I\otimes \tilde \Gamma\geq \varepsilon I  \, \}\nn}
 where $\epsilon$ is a constant taken sufficiently small. Thus, $\Cc_2$ is closed and bounded and thus compact. We conclude that $J_f$ admits minimum over $\Cc$. In a similar way, it is possible to prove that $J_g$ admits minimum over the set
 \al{\Cc=\{\, (\tilde \Pi,\tilde \Omega)\in \tilde{\Dc}   \hbox{ s.t. } \tilde \Pi\geq0,\, \tilde \Omega \geq0,\, \tilde \Pi\otimes I+I\otimes \tilde \Omega>0   \, \}\nn} with $\tilde{\Dc}=\Dc_{p_1}\times \Dc_{p_2m}$. 
We conclude that the dual problem admits solution. Accordingly, (\ref{opt_pI}) is the solution to the maximum entropy problem (\ref{ME_problem}). Finally, by defining $\Lambda=\tilde \Lambda^{-1}$, $\Gamma=\tilde \Gamma^{-1}$, $\Pi=\tilde \Pi^{-1}$ and $\Omega=\tilde \Omega^{-1}$ we obtain (\ref{optimal_K}). \qed
\end{proof}

We assumed that constraints (\ref{constr1}), (\ref{constr2}), (\ref{constr1f}) and (\ref{constr2f}) are totally binding in problem (\ref{ME_problem}). However, we are interested in the limiting cases where $c_{hj}^\prime$, $c_{kl}^{\prime\prime}$, $d_{h}^\prime$ and $d_{ki}^{\prime\prime}$ might be equal to zero in order to obtain a posterior corresponding to a Kronecker model. 
To include these scenarios, we consider the limits as $c_{hj}^\prime \rightarrow 0 $, $c_{kl}^{\prime\prime}\rightarrow 0$,
$d_{h}^\prime \rightarrow 0 $, $d_{ki}^{\prime\prime}\rightarrow 0$  and extend the maximum entropy solution by continuity. 

\pro  \label{extended_ME} Let $\Ec_1:=\{(h,j) \hbox{ s.t. } c_{hj}^\prime>0\}$, $\Ec_2:=\{(k,l) \hbox{ s.t. } c_{kl}^{\prime\prime}>0\}$.  If the sets $\Ec_1,\Ec_2$ are nonempty, then the maximum entropy solution extended by continuity is the probability density such that      
 $\check \theta_g$ is Gaussian, zero-mean and with kernel matrix as in (\ref{optimal_K}) where:
 \begin{itemize}
 \item $(i,j)\notin \Ec_1$ $\implies$ $\lambda_{hj}=0$, 
 \item $(k,l)\notin \Ec_2$ $\implies$ $\gamma_{kl}=0$.   
 \end{itemize}
\epro 
\begin{proof} First, recall that $\lambda_{hj}=\tilde \lambda_{hj}^{-1}$, $\gamma_{kl}=\tilde \gamma_{kl}^{-1}$ and $\tilde \lambda_{hj}$, $\tilde \gamma_{kl}$ are the original Lagrange multipliers. The maximum entropy solution $\check \theta_{ME}$ is such that, by Theorem \ref{teo_ME}, $\check g_{ME}^{[hk,jl]}$ is a zero-mean process with kernel matrix
\al{\frac{1}{\tilde \lambda_{hj}+\tilde\gamma_{kl}}P_{\Ic_g}.\nn} Moreover, if $\tilde \gamma_{kl}>0$ ($\lambda_{hj}>0$) then the maximum entropy solution $\check \theta_{ME}$ satisfies the corresponding constraint with equality. We consider two sequences $\{c_{hj}^{\prime(n)}\}_{n\in \Ns}$, $c_{hj}^{\prime (n)}>0$, and $\{c_{kl}^{\prime\prime(n)}\}_{n\in \Ns}$, $c_{kl}^{\prime\prime (n)}>0$, such that $c_{hj}^{\prime(n)}\rightarrow c_{hj}^{\prime}$ and $c_{kl}^{\prime\prime(n)}\rightarrow c_{kl}^{\prime\prime}$ as $n\rightarrow\infty$. For any $c_{hj}^{\prime(n)}$, we have
\al{\label{ineq1}c_{hj}^{\prime(n)}&\geq \sum_{k,l\in I_2}\Es[\|\check g_{ME}^{[hk,jl]}\|^2_{P_{\Ic_g}^{-1}}]\nn\\ &=\sum_{k,l\in I_2}\tr\left(\Es[\check g_{ME}^{[hk,jl]}\check g_{ME}^{[hk,jl]\top}]P^{-1}_{\Ic_g}\right)\nn\\
&= |\Ic_g| \sum_{k,l\in I_2}\frac{1}{\tilde \lambda_{hj}^{(n)}+\tilde \gamma_{kl}^{(n)}}>0.}In a similar way, we have 
\al{\label{ineq2}c_{kl}^{\prime\prime(n)}&\geq  |\Ic_g|  \sum_{h,j\in I_1}\frac{1}{\tilde \lambda_{hj}^{(n)}+\tilde \gamma_{kl}^{(n)}}>0.}
Let $(h,j)\notin \Ec_1$, then $c^{\prime(n)}_{hj}\rightarrow 0$. Then, we consider the arbitrary subsequences $\tilde \lambda^{(n_r)}_{hj}$ and $\tilde \gamma^{(n_r)}_{kl}$, with $r\in\Ns$, which admit limit. Taking into account (\ref{ineq1}), we may have two possible cases. {\em First case}. $\tilde\gamma_{kl}^{(n_r)}$ converges for some $k,l\in I_2$. By (\ref{ineq1}), this implies that $\tilde \lambda_{hj}^{(n_r)}\rightarrow \infty$ and thus $\lambda_{hj}^{(n_r)}\rightarrow 0$. {\em Second case}. $\tilde \gamma_{kl}^{(n_r)}\rightarrow \infty$ for any $k,l\in I_2$. This means that there exists $\bar r\in \Ns$ sufficiently large for which $\tilde \gamma_{kl}^{(n_r)}>0$ for any $k,l\in I_2$ and $r\geq \bar r$. Accordingly, the corresponding constraints in (\ref{ineq2}) are satisfied with equality for $r\geq \bar r$ and $n=n_r$. Hence $c_{kl}^{\prime \prime (n_r)}\rightarrow 0$ for any $k,l\in I_2$ and thus $\Ec_2=\emptyset$. This is not possible because $\Ec_2\neq \emptyset$ by assumption. We conclude that $c_{hj}^{\prime(n_r)}\rightarrow 0$ implies $\lambda_{hj}^{(n_r)}\rightarrow 0$ for any arbitrary subsequence which admits limit and thus $\lambda_{hj}^{(n)}\rightarrow 0$. Accordingly, in a similar way, it is possible to prove that $c_{kl}^{\prime\prime(n)}\rightarrow 0$ implies $\gamma_{kl}^{(n)}\rightarrow 0$. \qed
\end{proof}

\pro  \label{extended_ME2} Let $\Ac_1:=\{(h) \hbox{ s.t. } d_{h}^\prime>0\}$, $\Ac_2:=\{(k,i) \hbox{ s.t. } d_{ki}^{\prime\prime}>0\}$.  If the sets $\Ac_1,\Ac_2$ are nonempty, then the maximum entropy solution extended by continuity is the probability density such that      
$\check \theta_f$ is Gaussian, zero-mean and with kernel matrix as in (\ref{optimal_K}) where \begin{itemize} \item  $h\notin \Ac_1$ $\implies$ $\pi_{h}=0$,
\item $(k,l)\notin \Ac_2$ $\implies$ $\omega_{ki}=0$.   
\end{itemize}
\epro 

\begin{proof} The proof is similar to one of Proposition  \ref{extended_ME}. \qed
\end{proof}

Finally, in view of the Kolmogorov extension Theorem, from the probability density of $\check \theta$ we can characterize the probability law of $\theta$ maximizing the differential entropy.

\corr Consider the zero-mean Gaussian process $\theta=[\, \theta_g^\top\;\,\theta_f^\top]^\top$ with kernel, respectively, 
\al{ \label{ME_kernel}
\begin{aligned}   
  K_g&= X_g\otimes P\\
    K_f&= X_f\otimes R\\
  \end{aligned}
}
where $X_g$ and $X_f$ are defined as in (\ref{optimal_K}).   
For all finite sets $\Ic\subset \Ns\times \Ns$, its joint probability density is the extended solution to the maximum entropy problem (\ref{ME_problem}).\ecorr

It is worth noting that the maximum entropy solution is such that $g^{[hk,jl]}$ and $f^{[hk,i]}$ are zero-mean Gaussian processes with kernel function, respectively,
\al{\label{structK}K_g^{[hk,jl]}:=\frac{\lambda_{hj}\gamma_{kl}}{\lambda_{hj}+\gamma_{kl}}P, \; \; K_f^{[hk,i]}:=\frac{\pi_{h}\omega_{ki}}{\pi_{h}+\omega_{ki}}R.} In the case that both $\lambda_{hj}$ and $\gamma_{kl}$
 are equal to zero, we define by continuity $\lambda_{hj}\gamma_{kl}/(\lambda_{hj}+\gamma_{kl})=0$. Accordingly, $K_g^{[hk,jl]}=0$ if $\lambda_{hj}=0$ or $\gamma_{kl}=0$. Similarly, in the case that $\pi_{h}$ or $\omega_{ki}$
 is equal to zero, we define $\pi_{h}\omega_{ki}/(\pi_{h}+\omega_{ki})=0$.

In \cite{CHIUSO_PILLONETTO_SPARSE_2012} the transfer matrices in model (\ref{OEmodel_start}) are modeled as zero mean Gaussian processes with kernels 
\al{\label{structPillo}K_g^{[hk,jl]}:=\lambda_{hjkl}P, \; \; K_f^{[hk,i]}:=\pi_{hki} R}
where $\lambda_{hjkl},\pi_{hki}\geq 0$, with $h,j\in I_1$, $k,l\in I_2$ and $i\in I_u$, are the hyperparameters representing the so called scale factors. The prior in (\ref{structPillo}) has been used to learn sparse networks. Finally, it is worth noting that there are many more hyperparameters in (\ref{structPillo}) compared with (\ref{structK}).
\rema It is worth noting that the constants $c_{hj}^\prime$'s, $c_{kl}^{\prime\prime}$'s, $d_{h}^\prime$'s and $c_{ki}^{\prime\prime}$'s were needed only to derive the kernel structure in (\ref{ME_kernel}). Indeed, we are not interested in their particular values because, as we will show in the next section, the hyperparameters of the kernel will be found through another paradigm which exploits the measured data. \erema

{\em Hierarchical dynamic networks}. An interesting scenario included in this framework is the case of a system without input (i.e. $y(t)$ is a stationary stochastic process which is characterized by the transfer matrix $G(z)$) with $p_1=p_2$ and $E_1=E_2$. It is then natural to take \al{K_g=(\Lambda\otimes \Lambda)(\Lambda\otimes I_{p_1^2}+I_{p_1^2}\otimes \Lambda)^{-1}\otimes P\nn} and thus $g^{[hk,jl]}$ is modeled is a Gaussian process with kernel function
\al{\label{kernel_hier}K_g^{[hk,jl]}:=\frac{\lambda_{hj}\lambda_{kl}}{\lambda_{hj}+\lambda_{kl}}P.}  This model corresponds to a network which is  hierarchically organized into communities (clusters), the communities then grow, creating miniature copies of themselves. This idea can be applied recursively at the price that the dimension of $G(z)$ increases exponentially. Thus, $\Lambda$ is the Kronecker initiator for the prior on $G(z)$. These models have been proposed in  \cite{leskovec2010kronecker}: the main difference with respect to our approach is the fact that we do not impose the Kronecker decomposition on $G(z)$, but on $\mathrm{supp}(G(z))$; moreover $y(t)$ is not white noise (i.e. we consider a process with dynamics). 

\section{Hyperparameters Estimation}\label{sec:opt_procedure}

 In order  to compute $\hat \theta_g$ and $\hat \theta_f$, we have to estimate $\Sigma$ and $\zeta$. The latter contains the hyperparameters of $P, R$ and  $\xi:=\{\Lambda, \Gamma, \Pi,\Omega\}$.
The noise covariance matrix $\Sigma$ can be estimated  using a low-bias ARX-model as suggested in \citep{GOODWIN_1992}.
To estimate the hyperparameters of $P$ and $R$ we consider a nonparametric  unstructured model, i.e. we consider model (\ref{OEmodel_start}) but we do not require that $G(z)$ and $F(z)$ satisfy condition (\ref{cond_supp}).
Then, the hyperparameters of $P$
can be estimated by minimizing the negative log-marginal likelihood of $\yv^+$ computed using the unstructured model, see \cite{PILLONETTO_2011_PREDICTION_ERROR}. Throughout the paper, kernel $P$ is chosen as Stable Spline (SS) with a modulation factor \citep{ZORZI2018125}:{\small\al{ &[P]_{t,s}=\left[\frac{e^{-\beta(t+s)}e^{-\beta\max(t,s)}}{2}-\frac{e^{-3\beta\max(t,s)}}{6}\right]\cos(\omega_0 (t-s))
 \nn}}where $\beta\in(0,1)$ and $\omega_0\in (0,\pi)$ tune, respectively, the a priori information about the decay rate and the frequency content of the impulse responses. Kernel $R$ is chosen likewise. As explained in \cite{PILLONETTO_2011_PREDICTION_ERROR}, kernels of this type represent a good prior in the case that the impulse responses of the predictor have ``fast dynamics''.

Then, an estimate of $\xi$ is given by minimizing
the negative log-marginal likelihood $\ell$ of $\yv^+$ under model  (\ref{OEmodel_start}), (\ref{ME_kernel}) with $\Sigma$, $P$ and $R$ fixed as above. Under the assumptions in (\ref{approx_pb}), we have (up to constant terms) 
 \al{\label{loglik}
 \ell (\yv^+,\xi)&=\log \det V+\frac{1}{2}(\yv^+) ^\top V^{-1}\yv^+\nn\\
 V&= \Phi K_g\Phi^\top+\Psi K_f\Psi^\top+\Sigma \otimes I
}
where $K_g$ and $K_f$ have been defined in \eqref{ME_kernel} and $\xi$ belongs to the compact set 
\al{\label{def_XI}\Xi=\{(\Lambda,\Gamma,\Pi,\Omega)\in \Dc \hbox{ s.t. } 0\leq \Lambda,\Gamma,\Pi,\Omega\leq \kappa I\}}
where $\Dc=\Dc_{p_1^2} \times \Dc_{p_2^2} \times \Dc_{p_1} \times \Dc_{p_2 m}$ and $\kappa >0$ is a constant taken sufficiently large.  
Notice that the choice of $\Xi$ is not restrictive. Indeed, if we take: $\Gamma\geq 0$ such that $\gamma_{kl}>0$ for some $k,l\in I_2$; $\Lambda \geq 0$ such that $\lambda_{hj}\rightarrow \infty$ for some $h,j\in I_1$. Then, $\lambda_{hj}\gamma_{kl}/(\lambda_{hj}+\gamma_{kl})\rightarrow \infty$ and thus $\ell(\yv^+,\xi)\rightarrow \infty$. This means that the minimum cannot be attained for $\lambda_{hj}>\kappa$. In the case that $\Gamma=0$, then we have $\ell(\yv^+,\{\Lambda,\Gamma,\Pi,\Omega\})=\ell(\yv^+,\{\Lambda+\tilde \Lambda,\Gamma,\Pi,\Omega\})$ for any $\tilde\Lambda\in\Dc_{p_1^2}$ such that 
$\Lambda+\tilde \Lambda\geq 0$; accordingly, we can restrict the values of $\Lambda$ as in (\ref{def_XI}). A similar reasoning holds for $\Gamma$, $\Pi$ and $\Omega$. Since $\ell$ is a continuous function over the compact set $\Xi$, it follows that $\ell$ admits a point of minimum over $\Xi$.

The minimization of (\ref{loglik}) with respect to $\xi$ in the numerical examples presented in the next section is performed through the nonlinear optimization Matlab solver \textbf{fmincon.m} which uses a gradient based  method and the gradient, i.e. $\partial\; \ell (\mathrm y^+,\xi)\slash \partial \xi$,  is supplied. The computation of the gradient can be done in a similar way to that explained in \cite[Section 5]{CHEN20132213}.  In the case one is interested to make more efficient the algorithm, it is possible to perform the minimization by means of the scaled gradient projection method proposed in \cite{BONETTINI_2014}.  Finally, it is worth noting that only a local minimum can be found because $\ell(\mathrm y^+,\xi)$ is not convex with respect to $\xi$.

Next, we prove that the posterior mean of $\theta_g$ and $\theta_f$ under the prior defined through kernels in  (\ref{ME_kernel}) leads to a sparse Kronecker network. The analysis is performed by drawing inspiration from \cite{aravkin2014convex} which deals with the case of a sparse Bayesian network \citep{CHIUSO_PILLONETTO_SPARSE_2012}, i.e. the prior is defined through kernels in (\ref{structPillo}), and indeed our prior is closely related to Auto Relevance Determination (ARD) method \citep{mackay1994bayesian}. First, notice that the estimates $\hat \theta_g$ and $\hat \theta_f$ describe a sparse Kronecker network only if $\hat \xi=\{\hat \Lambda,\hat \Gamma, \hat \Pi,\hat \Omega\}$ is such that $\hat \Lambda,\Gamma, \hat \Pi,\hat \Omega$ are spare matrices. We assume that the noise variance at each node is the same, i.e. $\sigma_{hk}^2=\sigma^2$ with $h\in I_1$ and $k\in I_2$, and that 
\al{\label{cond_regressor}&(I\otimes P^{1/2})\Phi^\top\Phi(I\otimes P^{1/2}) =NI,\nn\\  &(I\otimes P^{1/2})\Phi^\top\Psi (I\otimes R^{1/2})=0,\nn\\ &(I\otimes R^{1/2})\Psi^\top\Psi (I\otimes R^{1/2})=NI.} Moreover  $\sigma^2$, $P$ are $R$ are fixed. First, we consider the case that $g^{[hk,jl]}$'s and $f^{[hk,i]}$'s are finite dimensional vectors, say $T$ their dimension. Accordingly, the corresponding transfer matrices $G(z)$ and $F(z)$ are finite impulse responses (FIR) and $T$ represents their practical length. Moreover, we assume that $N\geq T$. We define 
\al{\hat \theta_{g,LS}&=(\Phi^\top\Phi)^{-1}\Phi^\top \yv^+=\frac{1}{N}(I\otimes P)\Phi^\top \yv^+,\nn\\ \;\hat \theta_{f,LS}&=(\Psi^\top\Psi)^{-1}\Psi^\top \yv^+=\frac{1}{N}(I\otimes R)\Psi^\top \yv^+\nn}
which are the least squares estimators of $\theta_g$ and $\theta_f$. Thus, we have 
\al{ \label{LS_sparse}\hat g^{[hk,jl]}_{LS}=\frac{1}{N}P\phi^\top_{jl}\yv^+_{hk}, \; \; \hat f^{[hk,i]}_{LS}=\frac{1}{N}R\psi^\top_{i}\yv^+_{hk}.}

\pro \label{prop_ARD} Let $\lambda_{hj}$, $\gamma_{k,l}$, $\pi_h$ and $\omega_{ki}$ be the entries in the main diagonal of $\Lambda$, $ \Gamma$, $\Pi$, $ \Omega$, respectively. The following facts hold:
\begin{itemize}
\item if \al{\label{cond_ARD_1}\frac{T}{N}\sigma^2 -\|\hat g_{LS}^{[hk,jl]}\|^2_{  P^{-1}}\geq 0, \; \; \forall \, k,l\in I_2,} then the subset $$\Xi_1=\{(\Lambda,\Gamma,\Pi,\Omega) \hbox{ s.t. }\Lambda,\Gamma,\Pi,\Omega\geq 0, \; \lambda_{hj}=0\}$$ contains a local minimum for $\ell$; 
\item  if \al{\frac{T}{N}\sigma^2 -\|\hat g_{LS}^{[hk,jl]}\|^2_{  P^{-1}}\geq 0, \; \; \forall \, h,j\in I_1,\nn}then the subset $$\Xi_2=\{(\Lambda,\Gamma,\Pi,\Omega) \hbox{ s.t. }\Lambda,\Gamma,\Pi,\Omega\geq 0, \; \gamma_{kl}=0\}$$ contains a local minimum for $\ell$; 
\item  if \al{\frac{T}{N}\sigma^2 -\|\hat f_{LS}^{[hk,i]}\|^2_{ R^{-1}}\geq 0, \; \; \forall \,  k\in I_2,\;i\in I_u,\nn} then the subset $$\Xi_3=\{(\Lambda,\Gamma,\Pi,\Omega) \hbox{ s.t. }\Lambda,\Gamma,\Pi,\Omega\geq 0, \; \pi_{h}=0\}$$ contains a local minimum for $\ell$;
\item  if \al{\label{cond_ARD_fin}\frac{T}{N}\sigma^2 -\|\hat f_{LS}^{[hk,i]}\|^2_{  R^{-1}}\geq 0, \; \; \forall \,  h\in I_1,} then the subset $$\Xi_4=\{(\Lambda,\Gamma,\Pi,\Omega) \hbox{ s.t. }\Lambda,\Gamma,\Pi,\Omega\geq 0, \; \omega_{ki}=0\}$$ contains a local minimum for $\ell$.
\end{itemize}\epro
\begin{proof} Let 
\al{X:= \left[\begin{array}{cc}X_g\otimes I& 0 \\0 & X_f\otimes I\end{array}\right] , \; \;  \tilde K:=\left[\begin{array}{cc}I\otimes P & 0 \\0 & I\otimes R\end{array}\right],\nn} thus
\al{K:= \left[\begin{array}{cc}K_g& 0 \\0 & K_f\end{array}\right]=X  \tilde K=\tilde K  X= \tilde K^{1/2}X \tilde K^{1/2} \nn} and condition (\ref{cond_regressor}) can be written as 
\al{\label{cond_regressor2}  \tilde K^{1/2}  \left[\begin{array}{c}\Phi^\top \\ \Psi^\top\end{array}\right]\left[\begin{array}{cc}\Phi & \Psi \end{array}\right]  \tilde K^{1/2} =N I .}
In what follows, the symbol $\equiv$ means ``equal to up to terms not depending on $\xi$''. Next, we rewrite $\ell (\yv^+,\xi)$ using (\ref{LS_sparse}) and (\ref{cond_regressor2}).  First,
{\small \al{&\log\det(\Phi K_g \Phi^\top+\Psi K_f\Psi^\top+\sigma^2 I)\nn\\
&\equiv  \log\det\left(\frac{1}{\sigma^{2}}\left[\begin{array}{cc}\Phi & \Psi \end{array}\right]  K\left[\begin{array}{c}\Phi^\top \\ \Psi^\top\end{array}\right]+ I\right)\nn\\
&\equiv  \log\det\left(\frac{1}{\sigma^{2}}\left[\begin{array}{cc}\Phi & \Psi \end{array}\right]  \tilde K^{1/2}X^{1/2}X^{1/2}\tilde K^{1/2}\left[\begin{array}{c}\Phi^\top \\ \Psi^\top\end{array}\right]+ I\right)\nn\\
& \equiv \log\det\left(\frac{1}{\sigma^{2}} X^{1/2} \tilde K^{1/2}  \left[\begin{array}{c}\Phi^\top \\ \Psi^\top\end{array}\right]\left[\begin{array}{cc}\Phi & \Psi \end{array}\right]  \tilde K^{1/2} X^{1/2}+ I\right)\nn\\
 &= \log\det\left(\frac{N}{\sigma^{2}}\left[\begin{array}{cc}X_g\otimes I_T & 0 \\0 & X_f \otimes I_T\end{array}\right]   + I\right)\nn\\
 &=\log\det\left(\frac{N}{\sigma^{2}}X_g\otimes I_T+I\right)+\log\det\left(\frac{N}{\sigma^{2}}X_f\otimes I_T+I\right)\nn\\
 & = T\sum_{\substack{h\in I_1\\ k\in I_2}}\left[ \sum_{\substack{j\in I_1\\ l\in I_2}}\log \det \left(\frac{N}{\sigma^2}\frac{\lambda_{hj}\gamma_{kl}}{\lambda_{hj}+\gamma_{kl}}+1\right)\right.\nn\\ &
\left. \hspace{0.4cm}+\sum_{\substack{i\in I_u}}\log \det \left(\frac{N}{\sigma^2}\frac{\pi_{h}\omega_{ki}}{\pi_{h}+\omega_{ki}}+1\right) \right]\nn
}}
{\small \al{&(\yv^+)^\top(\Phi K_g \Phi^\top+\Psi K_f\Psi^\top+\sigma^2 I)^{-1}\yv^+ \nn\\
&=(\yv^+)^\top\left(\left[\begin{array}{cc}\Phi & \Psi \end{array}\right] K \left[\begin{array}{c}\Phi^\top \\ \Psi^\top\end{array}\right]+\sigma^2 I\right)^{-1}\yv^+ \nn\\
 &\equiv-\frac{1}{\sigma^{2}}(\yv^+)^\top\left[\begin{array}{cc}\Phi & \Psi \end{array}\right] \left(\left[\begin{array}{c}\Phi^\top \\ \Psi^\top\end{array}\right]\left[\begin{array}{cc}\Phi & \Psi \end{array}\right]+\sigma^2K^{-1} \right)^{-1} \left[\begin{array}{c}\Phi^\top \\ \Psi^\top\end{array}\right]\yv^+ \nn\\
 &=-\frac{N^2}{\sigma^{2}}\left[\begin{array}{cc}\hat \theta_{g,LS}^\top & \hat \theta_{f,LS}^\top \end{array}\right]  \tilde K^{-1/2}\left(\tilde K^{1/2} \left[\begin{array}{c}\Phi^\top \\ \Psi^\top\end{array}\right] \left[\begin{array}{cc}\Phi & \Psi \end{array}\right] \tilde K^{1/2} \right.\nn\\ &\left. \hspace{0.4cm}+\sigma^2 X^{-1}  \right)^{-1} \tilde K^{-1/2}\left[\begin{array}{c}\hat \theta_{g,LS} \\ \hat \theta_{f,LS}\end{array}\right]  \nn\\
  &=-\frac{N^2}{\sigma^{2}}\left[\begin{array}{cc}\hat \theta_{g,LS}^\top & \hat \theta_{f,LS}^\top \end{array}\right]   \tilde K^{-1/2} \left( NI    +\sigma^2 X^{-1}   \right)^{-1} \tilde K^{-1/2}\left[\begin{array}{c}\hat \theta_{g,LS} \\ \hat \theta_{f,LS}\end{array}\right]  \nn\\
    &=-\frac{N^2}{\sigma^{2}} \sum_{\substack{h\in I_1\\ k\in I_2}}\left[\sum_{\substack{j\in I_1\\ l\in I_2}} \frac{\|\hat g^{[hk,jl]}_{LS}\|_{P^{-1}}^2}{N+\sigma^2\frac{\lambda_{hj}+\gamma_{kl}}{\lambda_{hj}\gamma_{kl}}} +\sum_{i\in I_u}\frac{\|\hat f^{[hk,i]}_{LS}\|_{R^{-1}}^2}{N+\sigma^2\frac{\pi_{h}+\omega_{ki}}{\pi_{h}\omega_{ki}}} \right]\nn.
 }} Thus,
 {\small\al{ &\ell(\yv^+,\xi)\equiv \sum_{\substack{h\in I_1\\ k\in I_2}}\left[ T\sum_{\substack{j\in I_1\\ l\in I_2}}\log \det \left(\frac{N}{\sigma^2}\frac{\lambda_{hj}\gamma_{kl}}{\lambda_{hj}+\gamma_{kl}}+1\right)\right.\nn\\ &
  \hspace{0.4cm}+T\sum_{\substack{i\in I_u}}\log \det \left(\frac{N}{\sigma^2}\frac{\pi_{h}\omega_{ki}}{\pi_{h}+\omega_{ki}}+1\right)  \nn\\ &  \left.-\frac{N^2}{\sigma^{2}}   \sum_{\substack{j\in I_1\\ l\in I_2}} \frac{\|\hat g^{[hk,jl]}_{LS}\|_{P^{-1}}^2}{N+\sigma^2\frac{\lambda_{hj}+\gamma_{kl}}{\lambda_{hj}\gamma_{kl}}} -\frac{N^2}{\sigma^{2}}\sum_{i\in I_u}\frac{\|\hat f^{[hk,i]}_{LS}\|_{R^{-1}}^2}{N+\sigma^2\frac{\pi_{h}+\omega_{ki}}{\pi_{h}\omega_{ki}}} \right].\nn}}Next, we prove the first statement regarding $\lambda_{hj}$. The expression of the partial derivative of $\ell$ with respect to $\lambda_{hj}$ depends on the value of $\gamma_{kl}$'s. Given $\gamma_{kl}\geq 0$ for $k,l\in I_2$, let $\tilde J$ denotes the subset of $ I_2\times I_2$ containing the indexes $(k,l)$ such that  $\gamma_{kl}>0$. Then, we have:
{\small \al{& \frac{\partial \ell(\yv^+,\xi)}{\partial \lambda_{hj}}\nn\\ &=\sum_{\substack{(k, l)\in \tilde J}}\left(\frac{\gamma_{kl}}{\lambda_{hj}+\gamma_{kl}}\right)^2
\frac{T\left[\frac{\lambda_{hj}\gamma_{kl}}{\lambda_{hj}+\gamma_{kl}}+\frac{\sigma^2}{N}\right]-\| \hat g_{LS}^{[hk,jl]}\|^2_{P^{-1}}}{\left[\frac{\lambda_{hj}\gamma_{kl}}{\lambda_{hj}+\gamma_{kl}}+\frac{\sigma^2}{N}\right]^2},\nn
}} thus 
\al{\label{grad_ineq}\left.\frac{\partial \ell(\yv^+,\xi)}{\partial \lambda_{hj}}\right|_{\substack{\lambda_{hj}=0 }}= \sum_{\substack{(k, l)\in \tilde J}}
\frac{\frac{T }{N}\sigma^2-\| \hat g_{LS}^{[hk,jl]}\|^2_{P^{-1}}}{\sigma^4/N^2}.}
If condition (\ref{cond_ARD_1}) holds, and in view of (\ref{grad_ineq}), then \al{\label{cond_der}\left.\frac{\partial \ell(\yv^+,\xi)}{\partial \lambda_{hj}}\right |_{\lambda_{hj}=0}\geq 0.} Notice that $\ell$ is continuous in the compact set $\Xi_1$ and thus $\ell$ admits a point of minimum in $\Xi_1$.
Therefore, if (\ref{cond_der}) holds, then $\Xi_1$ contains a local minimum for $\ell$. The remaining statements can be proved in a similar way. \qed
\end{proof}

\rema \label{cor_ARD} In the nonparametric case, i.e. the limit case $T\rightarrow \infty$, the above result holds provided that the length of the data diverges in the same way. More precisely, assume that $T\rightarrow \infty$, $N\rightarrow \infty$ and $T/N\rightarrow 1$, then Proposition \ref{prop_ARD} holds where conditions (\ref{cond_ARD_1})-(\ref{cond_ARD_fin}) are replaced by 
\al{&\sigma^2 -\|\hat g_{LS}^{[hk,jl]}\|^2_{  P^{-1}}\geq 0, \; \; \forall \, k,l\in I_2\nn\\
&\sigma^2 -\|\hat g_{LS}^{[hk,jl]}\|^2_{  P^{-1}}\geq 0, \; \; \forall \, h,j\in I_1\nn\\
&\sigma^2 -\|\hat f_{LS}^{[hk,i]}\|^2_{  R^{-1}}\geq 0,  \; \; \forall \,  k\in I_2,\; i\in I_u\nn \\
&\sigma^2 -\|\hat f_{LS}^{[hk,i]}\|^2_{  R^{-1}}\geq 0, \; \; \forall \,  h\in I_1.\nn}\erema

Let $\hat G_{LS}(z)$ and $\hat F_{LS}(z)$ be the transfer matrices corresponding to the least squares estimates in (\ref{LS_sparse}). Then, the corresponding model is 
\al{  y(t)=\hat G_{LS}(z)y(t)+ \hat F_{LS}(z) +e(t)\nn}
The first statement of Remark \ref{cor_ARD} means that if the $\ell_2$ norms of $\hat g_{LS}^{[hk,jl]}$'s describing the conditional Granger causality relations from $y^\star_j$ to $y^\star_h$ do not 
dominate the variance of the noise process, then $\hat \lambda_{hj}=0$\footnote{Here, we implicitly assume that the local minimum in $\Xi_1$ is chosen as optimal solution.} and thus the posterior mean of $g^{[hk,jl]}$'s belonging to the block $(h,j)$ are null functions. Similar reasonings hold for the other statements.

\section{Numerical experiments}\label{section_simulation}
In this section we test the performance of the proposed estimator.  
\subsection{Synthetic examples}
We consider five Monte Carlo studies of 100 runs. For each run we generate an ARMAX model of the form \al{\label{model_true}y(t)=G_{\mathrm{TR}}(z)y(t)+F_{\mathrm{TR}}(z) u(t)+e(t)} as follows:
\begin{itemize}
\item the rational transfer matrices $G_{\mathrm{TR}}(z)$ and $F_{\mathrm{TR}}(z)$ are strictly casual and  generated through the MATLAB function \texttt{drmodel.m} where the order of the system is equal to $20$; moreover, their poles are restricted to have absolute value less than 0.95.
\item  some entries of $G_{\mathrm{TR}}(z)$ and $F_{\mathrm{TR}}(z)$ are set equal to the null function in such a way that they satisfy condition (\ref{cond_supp}); matrices  $E_1$, $E_2$, $A_1$ and $A_2$ are randomly generated in such a way such that the fraction of nonnull entries is equal to 0.6. 
\end{itemize}
For each run we generate a dataset of length $N=500$ if not otherwise specified; the input $u$ has uncorrelated components which are generated through the MATLAB function \texttt{idinput.m} as a realization from a random Gaussian noise with band [0, 0.4].

We compare the following estimators:
\begin{itemize}
\item \textbf{K}: this is the estimator based on model (\ref{OEmodel_start}) where the transfer matrices are modeled as zero mean Gaussian processes with kernels (\ref{ME_kernel});
\item \textbf{S}: this is the estimator based on model (\ref{OEmodel_start}) where the transfer matrices are modeled as zero mean Gaussian processes with kernels (\ref{structPillo});
\item \textbf{SS}: this is the estimator based on model (\ref{OEmodel_start}) where the transfer matrices are modeled as zero mean Gaussian processes with kernels $\lambda I_{p_1^2p_2^2}\otimes P$ and $\pi I_{p_1p_2m}\otimes R$ where the hyperparameters $\lambda,\pi\geq 0$ represent the scale factors. 
\end{itemize}

The impulse responses are truncated with length $T=50$, i.e. we make the approximation $G(z)= \sum_{t=1}^\infty G_t z^{-t} \approx \sum_{t=1}^T G_t z^{-t}$. Such a truncation does not affect the estimates because both $P$ and $R$ force the estimated impulse responses to decay exponentially, i.e. it does not involve any kind of trade-off between bias and variance \citep{PILLONETTO_DENICOLAO2010}. In practice, it is only important to take the practical length very large.     

The performance of these estimators is assessed through the following indexes:
\begin{itemize}
\item Average impulse response fit 
\al{\mathrm{AIRF}= \frac{fit(G_{\mathrm{TR}},\hat G)}{2}+\frac{fit(F_{\mathrm{TR}},\hat F)}{2} \nn}
where
\al{&fit(G_{\mathrm{TR}},\hat G)= 100\left(1-\sqrt{\frac{\frac{1}{T}\sum_{t=1}^T\|G_{\mathrm{TR},t}-\hat G_t\|^2}{\frac{1}{T}\sum_{t=1}^T\|G_{\mathrm{TR},t}-\bar G_{\mathrm{TR}}\|^2}} \right)\nn\\ 
&\bar G_{\mathrm{TR}}= \frac{1}{T}\sum_{t=1}^TG_{\mathrm{TR},t}\nn}
and $\hat G(z)$, $\hat F(z)$ are the estimated transfer matrices;
\item Fraction of misspecified edges{\small \al{\mathrm{ERR}=\frac{\|  E_1\otimes E_2 - \hat E_1\otimes \hat E_2\|_0}{2p_1^2p_2^2}+\frac{\|  A_1\otimes A_2 - \hat A_1\otimes \hat A_2\|_0}{2p_1p_2m}\nn}}
where $\mathrm{supp}(\hat G(z))=\hat E_1 \otimes \hat E_2$, $\mathrm{supp}(\hat F(z))=\hat A_1 \otimes \hat A_2$, while $E_1$, $E_2$, $A_1$ and $A_2$ correspond to the support of $G_{\mathrm{TR}}(z)$ and $F_{\mathrm{TR}}(z)$, respectively.
\end{itemize}

{\em First Monte Carlo study.} We consider the case in which model (\ref{model_true}) is without input $u$, i.e. we have a dynamic network corresponding to an ARMA stochastic process. In this case we only have to estimate $G(z)$, accordingly $\mathrm{AIRF}=fit(G_{\mathrm{TR}},\hat G)$ and $\mathrm{ERR}=\|  E_1\otimes E_2 - \hat E_1\otimes \hat E_2\|_0 /(p_1^2p_2^2)$. Besides the previous estimators, we also consider:
\begin{itemize}
\item \textbf{PEM}: this is the classical PEM approach that uses BIC for model order selection. The orders of the polynomials in the ARMA model ranges from 1 to 15 and are not allowed to be different from each other since this would lead to a combinatorial explosion of the number of competitive models.
\item \textbf{PEM+OR}: this is the same as \textbf{PEM} with an additional oracle knowing which impulse responses are null functions. 
\end{itemize} The length of the dataset is $N=3900$. It is worth noting the latter has been chosen so large in such a way that \textbf{PEM} with polynomial orders equal to 15 has enough data to estimate all the model parameters whose amount is 3848. 
\begin{figure}
\centering
\includegraphics[width=0.5\textwidth]{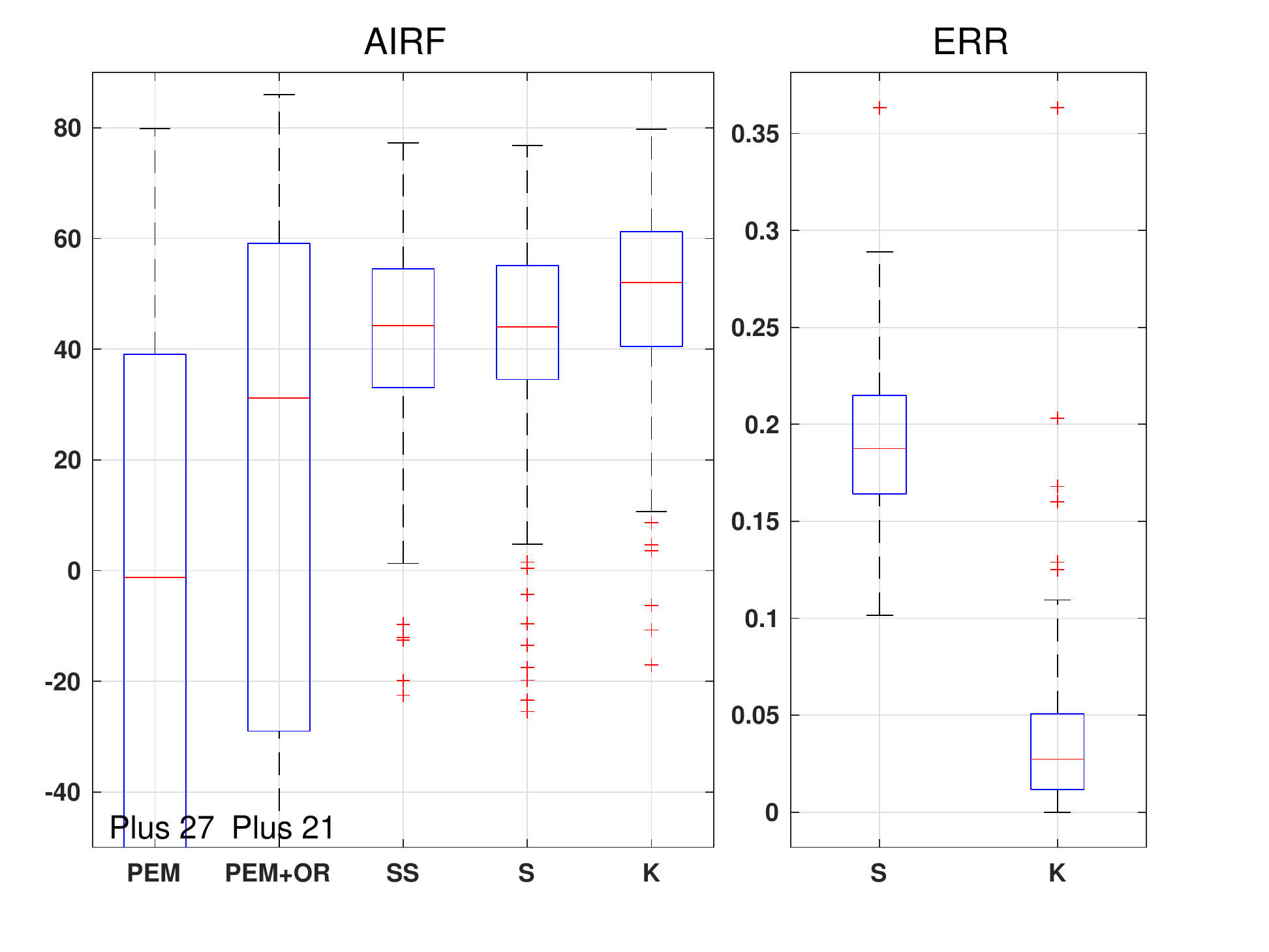}
\caption{First Monte Carlo study ``ARMA case'' with $p_1=4$, $p_2=4$.  {\em Left panel}. Average impulse response fit. {\em Right panel}. Fraction of misspecified edges.}\label{MC4}
\end{figure} Figure \ref{MC4} shows the performance of the estimators. \textbf{K} outperforms all the other estimators in terms of $\mathrm{AIRF}$ while both  \textbf{PEM+OR}, even though it owns the correct network topology, and \textbf{PEM} are the worst. Finally, \textbf{K} is better than \textbf{S} in terms of $\mathrm{ERR}$. Notice that we did not plot $\mathrm{ERR}$ for \textbf{PEM} and \textbf{SS} because they always provide a full network. In particular, in regard to \textbf{SS}, the entries of $G(z)$ share the same prior (and thus the same scale factor $\lambda$ whose estimated value is always strictly positive). We also considered Monte Carlo studies where $p_1$ and $p_2$ are different: the results, in terms of performance, do not change. 

{\em Second Monte Carlo study}. 
The output dimensions are $p_1=p_2=3$, while the input dimension is $m=2$. \begin{figure}
\centering
\includegraphics[width=0.5\textwidth]{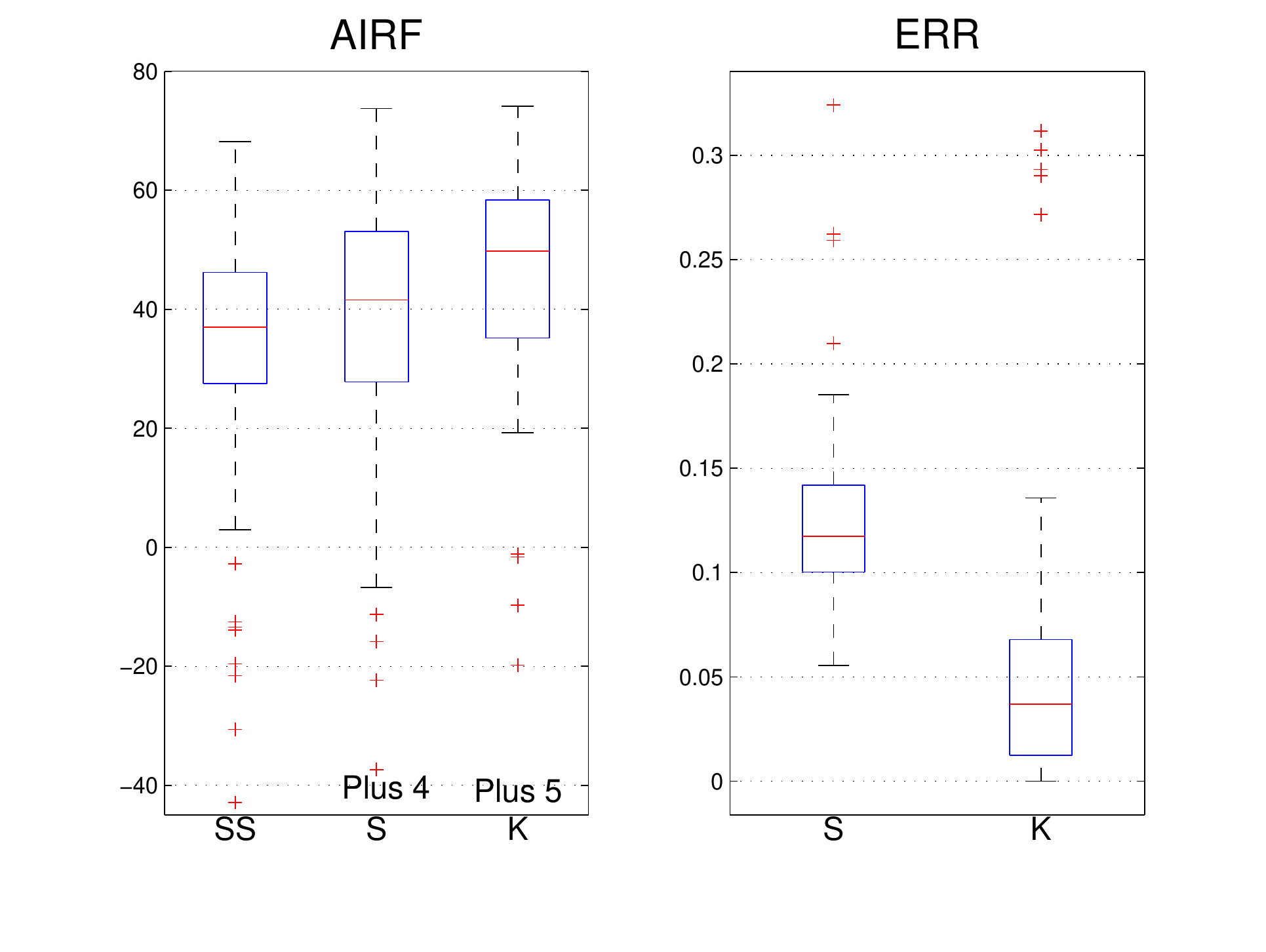}
\caption{Second Monte Carlo study with $p_1=3$, $p_2=3$ and $m=2$.  {\em Left panel}. Average impulse response fit. {\em Right panel}. Fraction of misspecified edges.}
\label{MC1}
\end{figure} Figure \ref{MC1} shows the performance of the three estimators: \textbf{K} outperforms all the other estimators in terms of $\mathrm{AIRF}$, while \textbf{SS} is the worst. Finally, \textbf{K} is better than \textbf{S} in terms of $\mathrm{ERR}$. 

{\em Third Monte Carlo study}.  We consider the case in which $p_1$ and $p_2$ are different; more precisely, we have $p_1=3$, $p_2=4$ and $m=2$.
\begin{figure}
\centering
\includegraphics[width=0.5\textwidth]{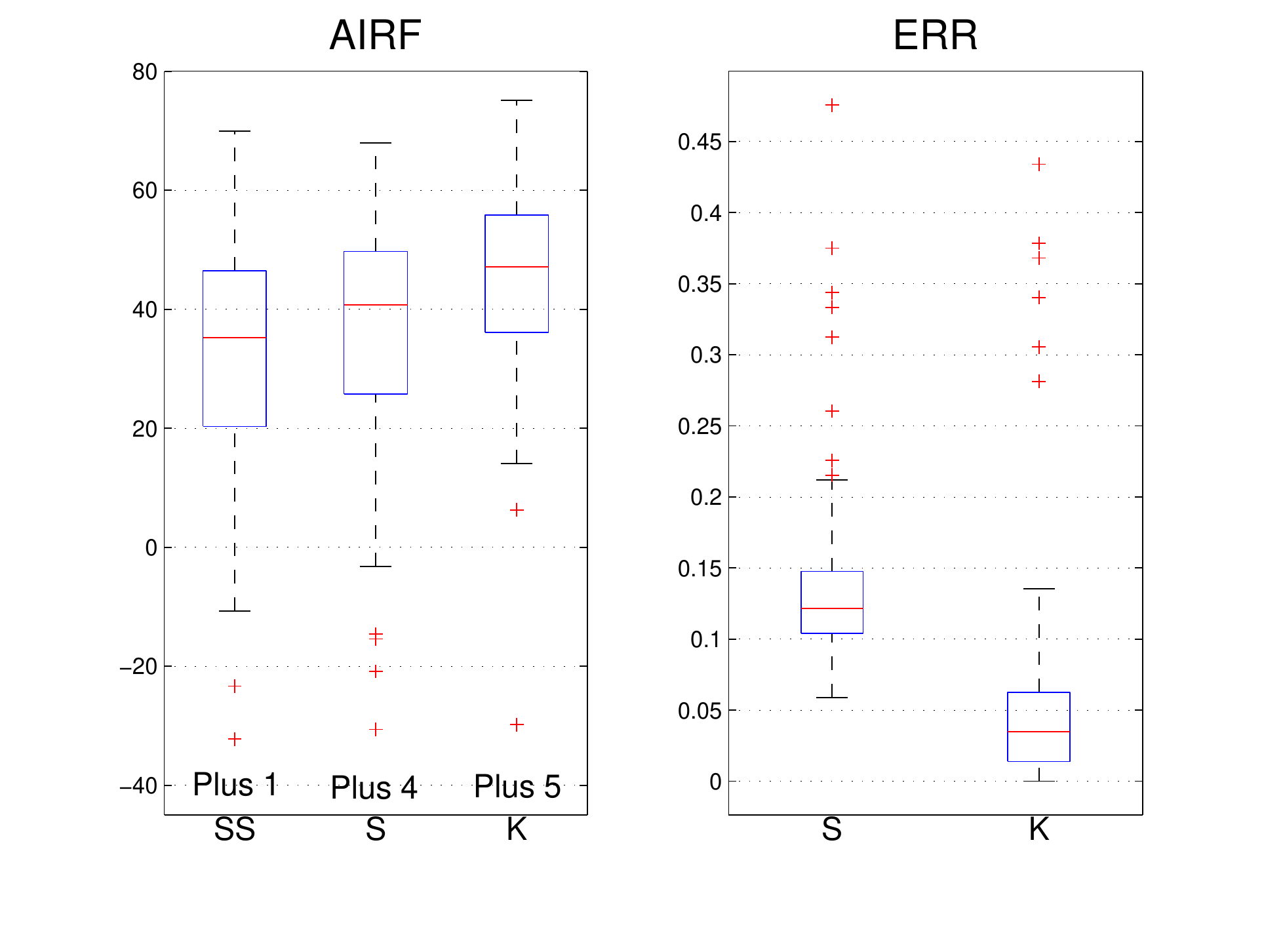}
\caption{Third Monte Carlo study with $p_1=3$, $p_2=4$ and $m=2$.  {\em Left panel}. Average impulse response fit. {\em Right panel}. Fraction of misspecified edges.}\label{MC2}
\end{figure}
Figure \ref{MC2} shows the obtained results. Also in this case \textbf{K} is the best estimator both in terms of $\mathrm{AIRF}$ and $\mathrm{ERR}$, while \textbf{SS} is the worst. 

{\em Fourth Monte Carlo study}. We consider $p_1=4$, $p_2=3$ and $m=2$; this Monte Carlo study differs from the last one because the output dimensions are swapped.
\begin{figure}
\centering
\includegraphics[width=0.5\textwidth]{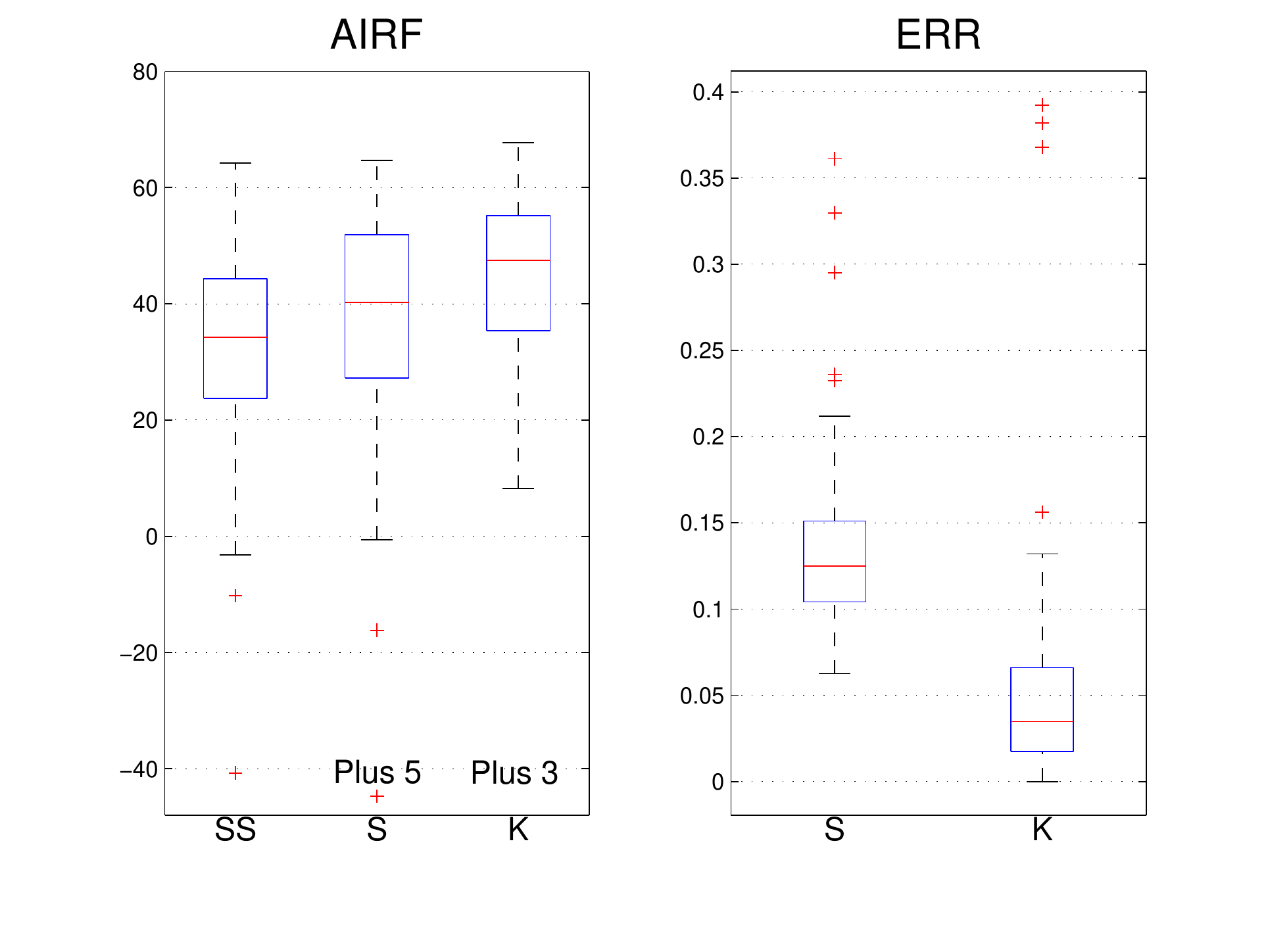}
\caption{Fourth Monte Carlo study with $p_1=4$, $p_2=3$ and $m=2$.  {\em Left panel}. Average impulse response fit. {\em Right panel}. Fraction of misspecified edges.}\label{MC3}
\end{figure} The situation substantially does not change, see Figure \ref{MC3}.

{\em Fifth Monte Carlo study}.  We consider the case in which model (\ref{model_true}) is without input $u$ and corresponds to a hierarchical dynamic network (i.e.  $E_1= E_2$). We  set $p_1=p_2=4$. Besides the previous estimators, we also consider:
\begin{itemize}
\item \textbf{H}: this is the estimator based on model (\ref{OEmodel_start}), without input, where the transfer matrix in $G(z)$ is modeled as zero mean Gaussian process with kernel (\ref{kernel_hier}).
\end{itemize} $\mathrm{AIRF}$ and $\mathrm{ERR}$ are defined as in the first Monte Carlo study.  \begin{figure}
\centering
\includegraphics[width=0.5\textwidth]{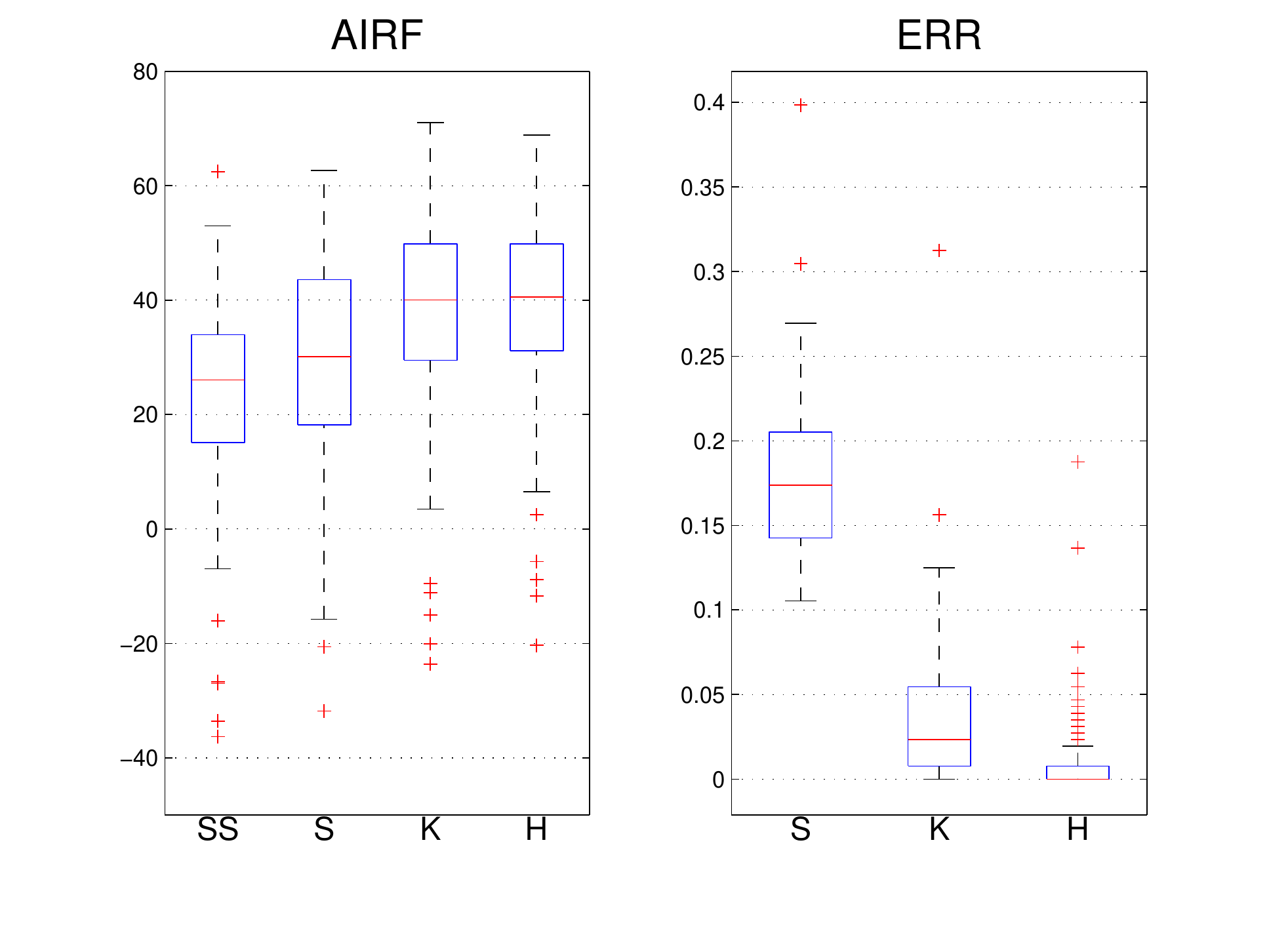}
\caption{Fifth Monte Carlo study ``ARMA case - hierarchical  network'' with $p_1=4$, $p_2=4$.   {\em Left panel}. Average impulse response fit. {\em Right panel}. Fraction of misspecified edges.}\label{MC5}
\end{figure} In Figure \ref{MC5} the performance of the estimators is depicted. \textbf{S} and \textbf{SS} are worse than \textbf{H} and \textbf{K}. The latter perform in a similar way in terms of $\mathrm{AIRF}$. \textbf{H}, however, outperforms \textbf{K} in terms of $\mathrm{ERR}$.  
 
\rema The negative log-marginal likelihood function is not convex with respect to the hyperparameters even in the simplest case where the hyperparameters are only the scale factors, see  \cite{6883125}. Accordingly, the kernel-based PEM methods
 could suffer about local minima in the estimation of the hyperparameters.  In particular, both \textbf{S} and \textbf{K}-\textbf{H} have the same issue. The main difference   is the number of hyperparameters to estimate: \textbf{K}-\textbf{H}  have to estimate very less hyperparameters than \textbf{S}. 
More precisely, in \textbf{K} we have to estimate $p_1^2+p_2^2+p_1+p_2m$ ($p_1^2+p_2^2$ in the case of systems without input) hyperparameters, in \textbf{H} we have to estimate $p_1^2$ hyperparameters (only applicable in the case of systems without input and $p_1=p_2$), while in 
\textbf{S} we have to estimate $p_1^2p_2^2+p_1p_2m$ ($p_1^2p_2^2$ in the case of systems without input) hyperparameters. Table \ref{tab:table1} summarizes the number of hyperparameters needed in the previous Monte Carlo studies. 
\begin{table}[h!]
  \begin{center}\caption{Number of hyperparameters to estimate.}
    \label{tab:table1}
   \begin{tabular}{l|c|c|c|} 
      \textbf{M.C. study} & \textbf{S} & \textbf{K}& \textbf{H}\\
      \hline
      First & 99 & 27 & -\\
      Second & 168 & 36& -\\
      Third & 168 & 35 &-\\
       Fourth   & 256 & 32 &- \\
Fifth        & 256 & 32 & 16 \\
    \end{tabular}
  \end{center} 
\end{table}
 Therefore, \textbf{K} and \textbf{H} alleviate the aforementioned issue since  
the search space for optimizing the marginal likelihood is small in respect to the one of \textbf{S}.\erema

\subsection{Learning a bike sharing system} \label{sec_bike} 
\begin{figure}
\centering
\includegraphics[width=0.2\textwidth]{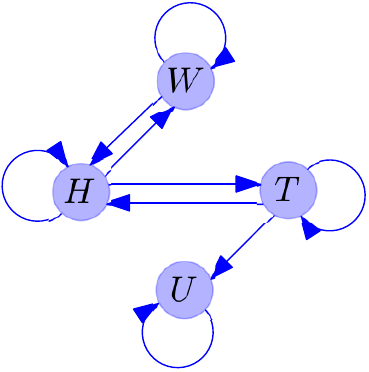}
\caption{Estimated network describing conditional Granger causality relations among: the number of users (U), wind speed (W), humidity (H) and apparent temperature (T).}\label{bike_E2}
\end{figure}

Bike sharing systems aim to make automatic bike rentals,  in particular   the  rental and the  return back. We consider the number of users and some variables describing the weather, that is wind speed, humidity and  apparent temperature. These quantities have been collected in the period 1 January 2011 - 14 May 2012 (500 days in total)  from the Capital Bike Sharing (CBS) system at Washington, D.C., USA. The sampling time is equal to 1 hour. For more details see \cite{fanaee2014event}. These data describe a four dimensional stochastic process $x(s)=[\,x_1(s) \; x_2(s) \; x_3(s) \; x_4(s) \, ]^\top$, where $x_1(s)$ is the number of users, $x_2(s)$ is the wind speed, $x_3(s)$ is the humidity and $x_4(s)$ is the apparent temperature. Process $x(s)$ is non-stationary during a day: for instance, we expect that the   number of users during the night  is smaller than the ones during the day. Accordingly, we consider the process
\al{y(t)=[\,x(24(t-1)+1)^\top \; \ldots  \; \,  x(24t)^\top \, ]^\top\nn}  taking values in $\Rs^{96}$ and the corresponding sampling time is 1 day. Thus, the corresponding dataset has length $N=500$. After detrending it, we model $y$ as the dynamic spatio-temporal model of Section \ref{sec_examples} where $p_1=24$, $p_2=4$. Then, we estimate $G(z)$ by \textbf{K} and the regression matrices are the ones in Remark \ref{remark_spatio}. The estimated network describing the conditional Granger causality relations among the number of users, wind speed, humidity and  apparent temperature is depicted in Figure \ref{bike_E2}. 
 As expected, the number of users does not conditionally Granger causes the weather variables. Moreover, the number of users is conditionally Granger caused by the apparent temperature.
\begin{figure}
\centering
\includegraphics[width=0.5\textwidth]{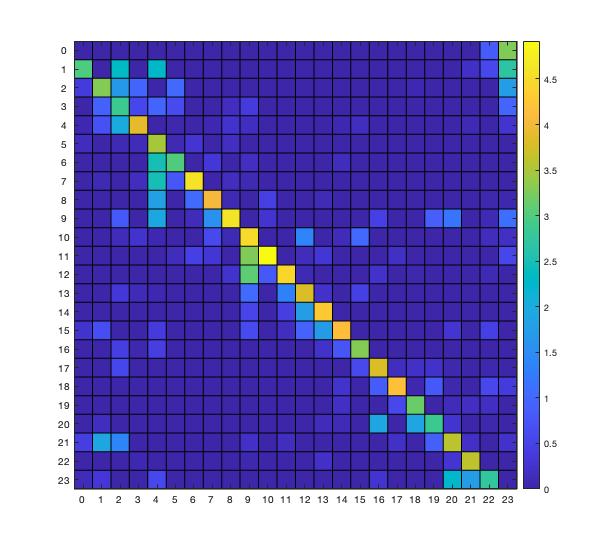}
\caption{Weighted matrix describing the estimated Granger causality relations among hours.}\label{fig_hours}
\end{figure} Figure \ref{fig_hours} shows the weighted matrix describing the conditional Granger causality relations among hours. More precisely, the entry in position $(h,j)$ is the $\ell_1$ norm of the $4\times 4$ impulse response, formed by $\hat g^{[hk,jl]}$ with $k,l\in \{1,2,3,4\}$, which describes the strength of the conditional Granger causality relation from hour $j$ to hour $h$. We can notice that all the hours are mainly conditionally Granger caused by the previous hour.

\section{Conclusions} \label{section_conclusions}
We have considered the problem to estimate a dynamic network describing conditional Granger causality relations and having a Kronecker structure. We have proposed a kernel-based PEM method for estimating such a model from data. The kernel functions have been derived from the maximum entropy principle. Numerical simulation showed that the proposed estimator outperforms the one estimating dynamic networks without structure. Finally, we have used the proposed method to learn the conditional Granger causality relations in a bike sharing system.

\bibliographystyle{model5-names}

\end{document}